\documentclass[12pt,leqno]{article}
\usepackage{amssymb}
\usepackage[mathscr]{eucal}
\usepackage{amsmath,amssymb,latexsym,theorem,bbm}
\newcommand{\EE}{\mathsf{E}}

\newcommand{\PP}{\mathsf{P}}

\newcommand{\RR}{\mathbb{R}}

\newcommand{\cB}{{\mathcal B}}

\newcommand{\cL}{{\mathcal L}}
\newcommand{\cN}{{\mathcal N}}

\newcommand{\dd}{\mathrm{d}}
\newcommand{\ee}{\mathrm{e}}

\DeclareMathOperator*{\argmax}{arg\,max}

\newcommand{\sign}{\operatorname{sign}}

\newcommand{\halpha}{\widehat{\alpha}}

\renewcommand{\leq}{\leqslant}
\renewcommand{\geq}{\geqslant}

\newcommand{\stoch}{\stackrel{\PP}{\longrightarrow}}
\newcommand{\distr}{\stackrel{\cL}{\longrightarrow}}
\newcommand{\distre}{\stackrel{\cL}{=}}

\newcommand{\bone}{\mathbbm{1}}

\newcommand{\proofend}{\hfill\mbox{$\Box$}}

\numberwithin{equation}{section}

\theoremstyle{change} \theorembodyfont{\em}
\newtheorem{Lem}{Lemma.}[section]
\newtheorem{Thm}[Lem]{Theorem.}

\theorembodyfont{\rm}
\newtheorem{Rem}[Lem]{Remark.}

\begin{document}

\begin{center}
 {\bfseries\Large Explicit formulas for Laplace transforms of certain} \\[2mm]
 {\bfseries\Large functionals of some time inhomogeneous diffusions}\\[5mm]

 {\sc\large M\'aty\'as $\text{Barczy}^{a,\diamond}$} {\large and}
 {\sc\large Gyula $\text{Pap}^{b}$}
\end{center}

\vskip0.2cm

$a$ University of Debrecen, Faculty of Informatics, Pf.~12, H--4010 Debrecen, Hungary

$b$ Bolyai Institute, University of Szeged, Aradi v\'ertan\'uk tere 1., H-6720 Szeged, Hungary

e--mails: barczy.matyas@inf.unideb.hu (M. Barczy), papgy@math.u-szeged.hu (G. Pap)

$\diamond$ Corresponding author

\vskip0.4cm




\renewcommand{\thefootnote}{}
\vspace*{0.2cm}
\footnote{The first author has been supported by the Hungarian
 Scientific Research Fund under Grant No.~OTKA T-079128/2009 
 and NKTH-OTKA-EU FP7 (Marie Curie action) co-funded 'MOBILITY' Grant No.~OMFB-00610/2010. 
The second author has been supported by the Hungarian Scientific Research Fund
 under Grant No.~OTKA T-079128/2009.}

\vspace*{-10mm}

\begin{abstract}
We consider a process \ $(X^{(\alpha)}_t)_{t\in[0,T)}$ \ given by the SDE
 \ $\dd X^{(\alpha)}_t = \alpha b(t) X^{(\alpha)}_t \, \dd t + \sigma(t) \, \dd B_t$,
 \ $t\in[0,T)$, \ with initial condition \ $X^{(\alpha)}_0=0$, \ where
 \ $T\in(0,\infty]$, \ $\alpha\in\RR$, \ $(B_t)_{t\in[0,T)}$ \ is a standard
 Wiener process, \ $b:[0,T)\to\RR\setminus\{0\}$ \ and \ $\sigma:[0,T)\to(0,\infty)$
 \ are continuously differentiable functions.
Assuming
 \ $
   \frac{\dd}{\dd t}\left(\frac{b(t)}{\sigma(t)^2}\right)
    =-2K\frac{b(t)^2}{\sigma(t)^2}
   $,
 \ $t \in[0,T)$, \ with some \ $K\in\RR$, \ we derive an explicit formula for the
 joint Laplace transform of
 \ $
    \int_0^t\frac{b(s)^2}{\sigma(s)^2}(X^{(\alpha)}_s)^2\,\dd s
   $
 \ and \ $(X^{(\alpha)}_t)^2$ \ for all \ $t\in[0,T)$ \ and for all \ $\alpha\in\RR$.
\ Our motivation is that the maximum likelihood estimator (MLE)  \ $\widehat\alpha_t$
 \ of \ $\alpha$ \ can be expressed in terms of these random variables.
As an application, we show that in case of \ $\alpha=K$, \ $K\ne0$,
 \[
   \sqrt{I_{K}(t)}\left(\widehat\alpha_t-K\right)
   \distre
   -\frac{\sign(K)}{\sqrt{2}}\frac{\int_0^1W_s\,\dd W_s}{\int_0^1(W_s)^2\,\dd s},
   \qquad \forall\;t\in(0,T) ,
 \]
 where \ $I_{K}(t)$ \ denotes the Fisher information for \ $\alpha$ \ contained
 in the observation \ $(X^{(K)}_s)_{s\in[0,\,t]}$, \ $(W_s)_{s\in[0,1]}$ \ is a standard Wiener
 process and \ $\distre$ \ denotes equality in distribution.
We also prove asymptotic normality of the MLE
 \ $\widehat\alpha_t$ \ of \ $\alpha$
 \ as \ $t\uparrow T$ \ for \ $\sign(\alpha - K) =  \sign(K)$, \ $K\ne0$.
\ As an example, for all \ $\alpha\in\RR$ \ and \ $T \in (0,\infty)$, \ we study
 the process \ $(X_t^{(\alpha)})_{t\in[0,T)}$ \ given by the SDE
 \ $\dd X_t^{(\alpha)}
    = - \frac{\alpha}{T-t} X_t^{(\alpha)} \, \dd t + \dd B_t$,
 \ $t\in[0,T)$, \ with initial condition \ $X_0^{(\alpha)}=0$.
\ In case of \ $\alpha>0$, \ this process is known as an $\alpha$-Wiener bridge,
 and in case of \ $\alpha=1$, \ this is the usual Wiener bridge.
\end{abstract}

\section{Introduction}

Several contributions have already been appeared containing explicit formulae
 for Laplace transforms of functionals of diffusion processes, see, e.g.,
 Borodin and Salminen \cite{BorSal}, Liptser and Shiryaev
 \cite[Sections 7.7 and 17.3]{LipShiII}, Arat\'o \cite{Ara2}, Yor \cite{Yor},
 Deheuvels and Martynov \cite{DehMar}, Deheuvels, Peccati and Yor \cite{DehPecYor},
 Mansuy \cite{Man}, Albanese and Lawi \cite{AlbLaw}, Kleptsyna and Le Breton \cite{KleLeb2},
 \cite{KleLeb3}, Hurd and Kuznetsov \cite{HurKuz}
 and Gao, Hannig, Lee and Torcaso \cite{GaoHanLeeTor}
 (the latter one is about the Laplace transform of the squared $L^2$-norm of some Gauss processes).
These formulae play an important role in theory of parameter estimation.
Most of the literature concern time homogeneous diffusion processes.

To describe our aims, let us start with the usual Ornstein--Uhlenbeck process
 \ $(Z_t^{(\alpha)})_{t\geq 0}$ \ given by the stochastic differential equation (SDE)
 \[
   \begin{cases}
    \dd Z_t^{(\alpha)}
    = \alpha \, Z_t^{(\alpha)} \, \dd t + \dd B_t ,
      \qquad t\geq0 ,\\
    \phantom{\dd} Z_0^{(\alpha)}=0 ,
  \end{cases}
 \]
where \ $\alpha\in\RR$ \ and \ $(B_t)_{t\geq 0}$ \ is a standard Wiener process.
An explicit formula is available for the Laplace transform of the random variable
 \ $\int_0^t(Z_s^{(\alpha)})^2\,\dd s$, \ $t\geq 0$,
 \ namely, for all \ $t\geq 0$ \ and \ $\mu>0$,
 \begin{align}\label{LP_OU_FORMULA}
   \EE\exp\left\{-\mu\int_0^t(Z_s^{(\alpha)})^2\,\dd s\right\}
   =\left( \frac{ \ee^{-\alpha t} \sqrt{\alpha^2+2\mu} }
                {\sqrt{\alpha^2+2\mu} \cosh(t\sqrt{\alpha^2+2\mu})
                 - \alpha \sinh(t\sqrt{\alpha^2+2\mu})} \right)^{\frac{1}{2}} ,
 \end{align}
 see, e.g., Liptser and Shiryaev \cite[Lemma 17.3]{LipShiII} or
 Gao, Hannig, Lee and Torcaso \cite[Theorem 4]{GaoHanLeeTor}.

Kleptsyna and Le Breton \cite[Proposition 3.2]{KleLeb2} presented an extension
 of the above mentioned result for fractional Ornstein--Uhlenbeck type processes.

In case of a time homogeneous diffusion process  \ $(H_t)_{t\geq 0}$,
 \ Albanese and Lawi \cite{AlbLaw} and Hurd and Kuznetsov \cite{HurKuz}
 recently addressed the question whether it is possible to compute the Laplace transform
 \[
    \EE\left[\ee^{-\int_0^t\phi(H_s)\,\dd s}q(H_t)\right], \qquad t>0,
 \]
 in an analitically closed form, where \ $\phi,q:\RR\to\RR$ \ are Borel measurable
 functions.
These papers provided a number of interesting cases when the Laplace transform
 can be evaluated in terms of special functions, such as hypergeometric
 functions.
Their methods are based on probabilistic arguments involving Girsanov theorem,
 and alternatively on partial differential equations involving Feynman--Kac
 formula.

As new results, in case of some time inhomogeneous diffusion processes,
 we will derive an explicit formula for the joint Laplace transform of certain functionals
 of these processes using the ideas of Florens-Landais and Pham \cite[Lemma 4.1]{FloPha},
 and see also Liptser and Shiryaev \cite[Lemma 17.3]{LipShiII}.
Let \ $T \in (0, \infty]$ \ be fixed.
Let \ $b: [0,T) \to \RR$ \ and \ $\sigma: [0,T) \to \RR$ \ be continuously
 differentiable functions.
Suppose that \ $\sigma(t) > 0$ \ for all \ $t\in[0,T)$, \ and \ $b(t)\ne0$ \ for
 all \ $t\in[0,T)$ \ (and hence \ $b(t)>0$ \ for all \ $t\in[0,T)$ \ or \ $b(t)<0$
 \ for all \ $t\in[0,T)$).
\ For all \ $\alpha \in \RR$, \ consider the process \ $(X_t^{(\alpha)})_{t\in[0,T)}$
 \ given by the SDE
 \begin{equation} \label{special_SDE}
  \begin{cases}
   \dd X_t^{(\alpha)}
    = \alpha \, b(t) X_t^{(\alpha)} \, \dd t + \sigma(t) \, \dd B_t ,
   \qquad t\in[0,T),\\
   \phantom{\dd} X_0^{(\alpha)}=0 .
  \end{cases}
 \end{equation}
The SDE \eqref{special_SDE} is a special case of Hull--White (or extended
 Vasicek) model, see, e.g., Bishwal \cite[page 3]{Bis}.
Assuming
 \begin{equation} \label{Laplace_feltetel}
   \frac{\dd}{\dd t}\left(\frac{b(t)}{\sigma(t)^2}\right)
   =-2K\frac{b(t)^2}{\sigma(t)^2}, \qquad t \in[0,T) ,
 \end{equation}
 with some \ $K\in\RR$, \ we derive an explicit formula for the joint Laplace
 transform of
 \begin{equation} \label{SEGED_abstract1}
    \int_0^t\frac{b(s)^2}{\sigma(s)^2}(X^{(\alpha)}_s)^2\,\dd s
    \qquad\text{and}\qquad
    (X^{(\alpha)}_t)^2
 \end{equation}
 for all \ $t\in[0,T)$ \ and for all \ $\alpha\in\RR$, \ see Theorem \ref{THEOREM_LAPLACE}.

We note that, using Lemma 11.6 in Liptser and Shiryaev \cite{LipShiII}, not assuming
 condition \eqref{Laplace_feltetel}, one can derive the following formula for the Laplace
 transform of
 \ $\int_0^t\frac{b(s)^2}{\sigma(s)^2}(X^{(\alpha)}_s)^2\,\dd s$,
 \begin{align*}
   \EE\exp\left\{-\mu\int_0^t\frac{b(s)^2}{\sigma(s)^2}(X^{(\alpha)}_s)^2\,\dd s\right\}
      = \exp\left\{\int_0^t\sigma(s)^2\gamma_t(s)\,\dd s\right\},
     \qquad  \mu>0,
 \end{align*}
 for all \ $t\in[0,T)$, \ where \ $\gamma_t:[0,t]\to \RR$ \ is the unique solution of
 the Riccati differential equation
 \begin{align}\label{LIPTSER_DE}
  \begin{cases}
   \frac{\dd \gamma_t}{\dd s}(s) = 2\mu \frac{b(s)^2}{\sigma(s)^2}-2\alpha b(s)\gamma_t(s)
                                  - \sigma(s)^2\gamma_t(s)^2,
        \qquad s\in[0,t],\\
   \phantom{\dd} \gamma_t(t)=0 .
  \end{cases}
 \end{align}
As a special case of our formula for the joint Laplace transform of \eqref{SEGED_abstract1},
 under the assumption \eqref{Laplace_feltetel}, we have an explicit formula for the Laplace transform
 of \ $\int_0^t\frac{b(s)^2}{\sigma(s)^2}(X^{(\alpha)}_s)^2\,\dd s$, \ $t\in[0,T)$,
 \ see Theorem \ref{THEOREM_LAPLACE} with \ $\nu=0$.
\ We suspect that, under the assumption \eqref{Laplace_feltetel}, the Riccati differential equation
 \eqref{LIPTSER_DE} may be solved explicitly.

We note that Deheuvels and Martynov \cite{DehMar} considered weighted Brownian
 motions \ $W_\gamma(t) := t^\gamma W_t$, \ $t \in (0,1]$, \ with
 \ $W_\gamma(0) := 0$, \ and weighted Brownian bridges
 \ $B_\gamma(t) := t^\gamma W_t - t^{\gamma+1} W_1$, \ $t \in (0,1]$,
 \ with \ $B_\gamma(t) := 0$, \ and \ with exponent \ $\gamma > -1$,
 \ where \ $(W_t)_{t\geq 0}$ \ is a standard Wiener process, and they
 explicitly calculated the Laplace transforms of the quadratic functionals
 \ $\int_0^1 W_\gamma(s)^2 \, \dd s$ \ and \ $\int_0^1 B_\gamma(s)^2 \, \dd s$
 \ by means of Karhunen--Lo\`eve expansions.
Deheuvels, Peccati and Yor \cite{DehPecYor} derived similar results for
 weighted Brownian sheets and bivariate weighted Brownian bridges.
Motivated by Theorems 1.3 and 1.4 in Deheuvels and Martynov \cite{DehMar} and Theorem 4.1
 in Deheuvels, Peccati and Yor \cite{DehPecYor}, we conjecture that
 our explicit formula in Theorem \ref{THEOREM_LAPLACE} for the joint Laplace transform
 of \eqref{SEGED_abstract1} may be expressed as an infinite product containing the
 eigenvalues of the integral operator associated with the covariance function
 of \ $(X^{(\alpha)}_t)_{t\in[0,T)}$.
\ Assumption \eqref{Laplace_feltetel} may play a crucial role in the calculation
 of these eigenvalues and also for deriving a (weighted) Karhunen--Lo\`eve expansion for
 \ $(X^{(\alpha)}_t)_{t\in[0,T)}$.
\ Once a (weighted) Karhunen--Lo\`eve expansion is available for \ $(X^{(\alpha)}_t)_{t\in[0,T)}$,
 \ one may derive the Laplace transform of
 \ $\int_0^t\frac{b(s)^2}{\sigma(s)^2}(X^{(\alpha)}_s)^2\,\dd s$, \ $t\in[0,T)$,
 \ as an infinite product.
We note that this approach can be carried through in the special case of a
 so-called $\alpha$-Wiener bridge with \ $\alpha=1/2$ \ (introduced and discussed later on).
Finally, we also remark that Gao, Hannig, Lee and Torcaso \cite{GaoHanLeeTor} used the same approach
 via Karhunen--Lo\`eve expansions for calculating the Laplace transform of the squared
 $L^2$-norm of some Gauss processes such as Ornstein-Uhlenbeck processes, time-changed
 Wiener bridges and integrated Wiener processes.

In Remark \ref{REM_technical} we give a third possible explanation for the role of the
 assumption \eqref{Laplace_feltetel}.

The random variables in \eqref{SEGED_abstract1} appear in the maximum
 likelihood estimator (MLE) \ $\widehat\alpha_t$ \  of \ $\alpha$ \ based on
 an observation \ $(X^{(\alpha)}_s)_{s\in[0,\,t]}$.
\ This is the reason why it is useful to calculate their joint Laplace
 transform explicitly.
For a more detailed discussion, see Sections \ref{Section_Maximum_Likelihood}
 and \ref{Section_Example}.

It is known that, under some conditions on \ $b$ \ and \ $\sigma$ \ (but
 without assumption \eqref{Laplace_feltetel}), the distribution of the MLE
 \ $\widehat\alpha_t$ \ of \ $\alpha$ \ normalized by
 Fisher information can converge to the standard normal distribution, to the
 Cauchy distribution or to the distribution of
 \ $c \int_0^1 W_s \, \dd W_s \, \big/ \int_0^1(W_s)^2 \, \dd s$, \ where
 \ $(W_s)_{s\in[0,1]}$ \ is a standard Wiener process, and \ $c=1/\sqrt{2}$ \ or
 \ $c=-1/\sqrt{2}$, \ see Luschgy \cite[Section 4.2]{Lus1} and Barczy and
 Pap \cite{BarPap}.
As an application of the joint Laplace transform of \eqref{SEGED_abstract1},
 under the conditions \ $\int_0^T\sigma(s)^2\,\dd s<\infty$ \ and
 \begin{align}\label{SEGED_abstract2}
  b(t)=\frac{\sigma(t)^2}{-2K\int_t^T\sigma(s)^2\,\dd s},
        \qquad t\in[0,T),
 \end{align}
 with some \ $K\ne0$ \ (note that in this case condition \eqref{Laplace_feltetel} is satisfied),
 we give an alternative proof for
 \[
   \sqrt{I_{\alpha}(t)}
   \Big(\widehat\alpha_t-\alpha\Big)
   \distr
   \begin{cases}
    \cN(0,1)  & \text{if \ $\sign(\alpha - K) = \sign(K)$,}\\[1mm]
    -\frac{\sign(K)}{\sqrt{2}} \,
    \frac{\int_0^1 W_s \, \dd W_s}{\int_0^1 (W_s)^2 \, \dd s}
              & \text{if \ $\alpha=K$,}
    \end{cases}
   \qquad\text{as \ $t\uparrow T,$}
 \]
 where \ $I_{\alpha}(t)$ \ denotes the Fisher information for \ $\alpha$
 \ contained in the observation \ $(X^{(\alpha)}_s)_{s\in[0,\,t]}$, \ $(W_s)_{s\in[0,1]}$
 \ is a standard Wiener process and \ $\distr$ \ denotes convergence in distribution,
 see Theorem \ref{THM_alpha_1/2_general}.
\ In fact, in case of \ $\alpha=K$, \ for all \ $t\in(0,T)$,
 \[
   \sqrt{I_{K}(t)}\left(\widehat\alpha_t-K\right)
   \distre
   -\frac{\sign(K)}{2\sqrt{2}}\frac{(W_1)^2-1}{\int_0^1(W_s)^2\,\dd s}
   = -\frac{\sign(K)}{\sqrt{2}}\frac{\int_0^1W_s\,\dd W_s}{\int_0^1(W_s)^2\,\dd s},
 \]
 where \ $\distre$ \ denotes equality in distribution, see Theorem \ref{THM_alpha_1/2_general}.
We note that in case of \ $\sign(\alpha - K) = - \sign(K)$, \ one can prove
\ $\sqrt{I_{\alpha}(t)}
   \left(\widehat\alpha_t-\alpha\right)
   \distr \zeta$
 \ as \ $t\uparrow T$, \ where \ $\zeta$ \ is a random variable with standard Cauchy distribution,
 see, e.g., Luschgy \cite[Section 4.2]{Lus1} or Barczy and Pap \cite{BarPap}.
The proof in this case is based on a martingale limit theorem, and we do not
 know whether one can find a proof using the explicit form of the joint Laplace
 transform of \eqref{SEGED_abstract1}.

By Barczy and Pap \cite[Corollaries 9 and 11]{BarPap}, under the conditions
 \ $\int_0^T\sigma(s)^2\,\dd s<\infty$ \ and \eqref{SEGED_abstract2}, we have
 for all \ $\alpha\ne K$, \ $K\ne0$, \ the MLE \ $\widehat\alpha_t$ \ of \ $\alpha$ \ is
 asymptotically normal with an appropriate \emph{random} normalizing factor,
 see also Remark \ref{REM_BARCZY_PAP}.
In case of \ $\alpha=K$, \ $K\ne0$, \ under the above conditions,
 we determine the distribution of this randomly normalized MLE using the joint
 Laplace transform of \eqref{SEGED_abstract1}, see Theorem \ref{THM_alpha_1/2_general_random}.
As a by-product of this result, giving a counterexample, we show that Remark
 1.47 in Prakasa Rao \cite{Rao} contains a mistake, see Remark \ref{REM_RAO}.

Using the explicit form of the Laplace transform we also prove strong consistency of
 the MLE of \ $\alpha$ \ for all \ $\alpha\in\RR$, \ see Theorem \ref{THM_CONSISTENCY}.

As an example, for all \ $\alpha\in\RR$ \ and \ $T \in (0,\infty)$, \ we study
 the process \ $(X_t^{(\alpha)})_{t\in[0,T)}$ \ given by the SDE
 \begin{equation} \label{alpha_W_bridge}
   \begin{cases}
   \dd X_t^{(\alpha)}
    = - \frac{\alpha}{T-t} X_t^{(\alpha)} \, \dd t + \dd B_t ,
      \qquad t\in[0,T),\\
    \phantom{\dd} X_0^{(\alpha)}=0 .
  \end{cases}
 \end{equation}
In case of \ $\alpha>0$, \ this process is known as an $\alpha$-Wiener bridge,
 and in case of \ $\alpha=1$, \ this is the usual Wiener bridge.
As a special case of the explicit form of the joint Laplace transform
 of \eqref{SEGED_abstract1}, we obtain the joint Laplace transform of
 \ $\int_0^t\frac{(X_u^{(\alpha)})^2}{(T-u)^2}\,\dd u$ \ and \ $(X_t^{(\alpha)})^2$ \
 for all \ $t\in[0,T)$, \ see Theorem \ref{Psi}.
As a special case of this latter formula we get the Laplace transform of
 \ $\int_0^t\frac{(B_u)^2}{(T-u)^2}\,\dd u$, \ $t\in[0,T)$, \ which was first
 calculated by Mansuy \cite[Proposition 5]{Man}, see Remark \ref{REM_MANSUY}.
Finally, we remark that in case of \ $\alpha>0$ \ unweighted and weighted
 Karhunen--Lo\`eve expansions are available for the \ $\alpha$-Wiener bridge
 \ $(X^{(\alpha)}_t)_{t\in[0,T)}$ \ on \ $[0,T]$ \ and \ $[0,S]$ \ with \ $0<S<T$,
 \ respectively, see Barczy and Igl\'oi \cite{BarIgl}.
Further, using the weighted Karhunen--Lo\`eve expansion, one can also get the Laplace
 transform of \ $\int_0^t\frac{(X^{(1/2)}_s)^2}{(T-s)^2}\,\dd s$, $t\in[0,T)$,
 \ see Barczy and Igl\'oi \cite[Proposition 3.1]{BarIgl}, i.e., in the special case
 of an \ $\alpha$-Wiener bridge with \ $\alpha=1/2$ \ the approach using Karhunen--Lo\`eve expansions
 mentioned earlier can be carried through.

\section{Laplace transform}\label{Section_Laplace}

Let \ $T \in (0, \infty]$ \ be fixed.
Let \ $b: [0,T) \to \RR$ \ and \ $\sigma: [0,T) \to \RR$ \ be continuously
 differentiable functions.
Suppose that \ $\sigma(t) > 0$ \ for all \ $t\in[0,T)$, \ and \ $b(t)\ne0$ \
 for all \ $t\in[0,T)$ \ (and hence \ $b(t)>0$ \ for all \ $t\in[0,T)$ \ or
 \ $b(t)<0$ \ for all \ $t\in[0,T)$).
\ For all \ $\alpha \in \RR$, \ consider the SDE \eqref{special_SDE}.
Note that the drift and diffusion coefficients of the SDE \eqref{special_SDE}
 satisfy the local Lipschitz condition and the linear growth condition
 (see, e.g., Jacod and Shiryaev \cite[Theorem 2.32, Chapter III]{JacShi}).
By Jacod and Shiryaev \cite[Theorem 2.32, Chapter III]{JacShi}, the SDE
 \eqref{special_SDE} has a unique strong solution
 \begin{align}\label{SolX}
   X^{(\alpha)}_t
   = \int_0^t
      \sigma(s) \exp\left\{ \alpha \int_s^t b(u) \, \dd u \right\} \dd B_s ,
   \qquad t \in [0,T) .
 \end{align}
Note that \ $(X_t^{(\alpha)})_{t \in [0,T)}$ \ has continuous sample paths by the
 definition of strong solution, see, e.g., Jacod and Shiryaev
 \cite[Definition 2.24, Chapter III]{JacShi}.
For all \ $\alpha \in \RR$ \ and \ $t \in (0,T)$, \ let \ $\PP_{X^{(\alpha)},\,t}$
 \ denote the distribution of the process \ $(X_s^{(\alpha)})_{s \in [0,\,t]}$ \ on
 \ $\big(C([0,t]),\cB(C([0,t]))\big)$, \ where \ $C([0,t])$ \ and \ $\cB(C([0,t]))$ \ denote
 the set of all continuous real valued functions defined on \ $[0,t]$ \ and the Borel
 $\sigma$-field on \ $C([0,t])$, \ respectively.
The measures \ $\PP_{X^{(\alpha)},\,t}$ \ and  \ $\PP_{X^{(\beta)},\,t}$ \ are equivalent
 for all \ $\alpha$, \ $\beta\in\RR$ \ and for all \ $t\in(0,T)$, \ and
 \begin{align}\label{SEGED_Radon_Nicodym}
   \frac{\dd \PP_{X^{(\alpha)},\,t}}{\dd \PP_{X^{(\beta)},\,t}}
       \left(X^{(\beta)}\big\vert_{[0,t]}\right)
       =\exp \left\{ (\alpha-\beta)
                      \int_0^t
                        \frac{b(s)}{\sigma(s)^2} X^{(\beta)}_s \, \dd X_s^{(\beta)}
                       -\frac{\alpha^2-\beta^2}{2}
                        \int_0^t
                         \frac{b(s)^2}{\sigma(s)^2}(X^{(\beta)}_s)^2
                          \, \dd s \right\},
 \end{align}
 see, e.g., Liptser and Shiryaev \cite[Theorem 7.19]{LipShiI}.
Note also that for all \ $s\in[0,T)$, \ $X_s^{(\alpha)}$ \ is normally distributed with
 mean \ $0$ \ and with variance
 \begin{align}\label{SEGED17_uj}
    V(s;\alpha)
     := \EE \big( X_s^{(\alpha)} \big)^2
      = \int_0^s
         \sigma(u)^2
         \exp\left\{ 2 \alpha \int_u^s b(v) \, \dd v \right\} \, \dd u ,
      \qquad s \in [0,T) ,
 \end{align}
 and then, by the conditions on \ $b$ \ and \ $\sigma$, \ $V(s;\alpha)>0$
 \ for all \ $s\in(0,T)$.

The next lemma is about the solutions of the differential equation (DE) \eqref{Laplace_feltetel}.

\begin{Lem}\label{LEMMA_DE}
Let \ $T\in(0,\infty]$ \ be fixed and let \ $b:[0,T)\to\RR\setminus\{0\}$
 \ and \ $\sigma:[0,T)\to(0,\infty)$ \ be continuously differentiable functions.
The DE \eqref{Laplace_feltetel} leads to a Bernoulli type DE
 having solutions
 \begin{align}\label{Laplace_feltetel2}
   b(t)=\frac{\sigma(t)^2}{2\left(K\int_0^t\sigma(s)^2\,\dd s+C\right)},
       \qquad t\in[0,T),
 \end{align}
 where \ $C\in\RR$ \ is such that the denominator \ $K\int_0^t\sigma(s)^2\,\dd s+C\ne0$ \
 for all \ $t\in[0,T)$.
\end{Lem}

\noindent{\bf Proof.}
The DE \eqref{Laplace_feltetel} can be written in the form
 \begin{align*}
     \frac{b'(t)\sigma(t)-2b(t)\sigma'(t)}{\sigma(t)^3}
        =-2K\frac{b(t)^2}{\sigma(t)^2}, \qquad t\in[0,T),
 \end{align*}
 which is equivalent to the Bernoulli type DE
 \begin{align*}
    b'(t)-2b(t)(\ln(\sigma(t)))'=-2K b(t)^2,\qquad t\in[0,T).
 \end{align*}
Since \ $b(t)\ne 0$ \ for all \ $t\in[0,T)$, \ we get
 \[
   b'(t)b(t)^{-2}-2(\ln(\sigma(t)))'b(t)^{-1}
      =-2K,\qquad t\in[0,T).
 \]
Let \ $u(t):=b(t)^{-1}$, \ $t\in[0,T)$. \ Then
 \begin{align}\label{DE_HOMOGENEOUS}
    -u'(t)-2(\ln(\sigma(t)))'u(t)=-2K,\qquad t\in[0,T),
 \end{align}
 which is an inhomogeneous linear differential equation.
The homogeneous linear DE \ $v'(t)+2(\ln(\sigma(t)))'v(t)=0$ \
 has solutions \ $v(t)=2C\sigma(t)^{-2}$, \ $t\in[0,T)$, \ $C\in\RR$, \
 and hence
 \[
   u(t)=2K\frac{\int_0^t\sigma(s)^2\,\dd s}{\sigma(t)^2},\qquad t\in[0,T),
 \]
 is a  particular solution of the inhomogeneous linear DE \eqref{DE_HOMOGENEOUS},
 which yields the assertion.
\proofend

Now we derive an explicit formula for the joint Laplace transform
 of \ $\int_0^t\frac{b(s)^2}{\sigma(s)^2}\big( X^{(\alpha)}_s \big)^2 \, \dd s$
 \ and \ $\big(X^{(\alpha)}_t\big)^2$
 \ for all \ $t\in[0,T)$ \ under the assumption \eqref{Laplace_feltetel2} on \ $b$ \ and \ $\sigma$.
\ We use the same technique (sometimes called Novikov's method, see, e.g., Arat\'o \cite{Ara2})
 as in the proof of Lemma 4.1 in Florens-Landais and Pham \cite{FloPha} or see also
 the proof of Lemma 17.3 in Liptser and Shiryaev \cite{LipShiII}.

\begin{Thm}\label{THEOREM_LAPLACE}
Let \ $(X_t^{(\alpha)})_{t\in[0,T)}$ \ be the process given by the SDE \eqref{special_SDE}
 where \ $b$ \ is given by \eqref{Laplace_feltetel2}.
Then for all \ $\mu>0$, \ $\nu\geq 0$, \ and \ $t\in[0,T)$, \ we have
 \begin{align*}
   &\EE\exp\left\{ - \mu \int_0^t \frac{b(s)^2}{\sigma(s)^2}(X^{(\alpha)}_s)^2 \, \dd s
                   -\nu[X_t^{(\alpha)}]^2
           \right\} \\[2mm]
   &=
          \frac{B_{K,C}(t)^{\frac{K-\alpha}{4}}}
            {\sqrt{\cosh\!\left(\!\frac{\sqrt{2\mu+(\alpha-K)^2}}{2}\ln(B_{K,C}(t))
                  \!\right)
            \! -\frac{\alpha-K-4\nu\left(K\int_0^t\sigma(s)^2\,\dd s+C\right)}
               {\sqrt{2\mu+(\alpha-K)^2}}
              \sinh\!\left(\!\frac{\sqrt{2\mu+(\alpha-K)^2}}{2}\ln(B_{K,C}(t))
                  \!\right)}},
 \end{align*}
 where
 \begin{align*}
   B_{K,C}(t)
    :=\begin{cases}
        \left(1+\frac{K}{C}\int_0^t\sigma(s)^2\,\dd s\right)^{\frac{1}{K}}
           & \quad \text{if \ $K\ne 0$,}\\[2mm]
        \exp\left\{\frac{1}{C}\int_0^t\sigma(s)^2\,\dd s\right\}
           &  \quad \text{if \ $K=0$,}
     \end{cases}
     \qquad t\in[0,T).
 \end{align*}
\end{Thm}

For the proof of Theorem \ref{THEOREM_LAPLACE} we need two lemmas.
The first one can be considered as a preliminary version of Theorem \ref{THEOREM_LAPLACE},
 the second one is about the variance of \ $X_t^{(\alpha)}$.

\begin{Lem}\label{Psi_general}
Let \ $(X_t^{(\alpha)})_{t\in[0,T)}$ \ be the process given by the SDE \eqref{special_SDE}.
If assumption \eqref{Laplace_feltetel} is satisfied with some \ $K\in\RR$ \
 and if \ $\sign(b)=\pm\bone_{[0,T)}$, \ then for all \ $\mu>0$, \ $\nu\geq 0$
 \ and \ $t \in [0,T)$, \ we have
 \begin{align}\label{LAPLACE}
   \EE \exp \left\{ - \mu \int_0^t \frac{b(s)^2}{\sigma(s)^2}(X^{(\alpha)}_s)^2 \, \dd s
                   -\nu[X_t^{(\alpha)}]^2
           \right\}
    =
     \left(
     \frac{\exp\left\{-A_{\mu,\alpha,K}^{\pm}\int_0^tb(s)\,\dd s\right\}}
         {1+\left(2\nu-A_{\mu,\alpha,K}^{\pm}\frac{b(t)}{\sigma(t)^2}\right)
               V(t;\alpha-A_{\mu,\alpha,K}^{\pm})}
    \right)^{\frac{1}{2}},
 \end{align}
 where \ $A_{\mu,\alpha,K}^{\pm}:=\alpha-K\mp\sqrt{2\mu+(\alpha-K)^2}$.
\end{Lem}

\noindent{\bf Proof.}
For all \ $\mu>0$, \ $\nu\geq 0$ \ and \ $t\in[0,T)$, \ let
 \[
   \Psi_t(\alpha,\mu,\nu)
     :=\EE\left(
        \exp
         \left\{
           - \mu \int_0^t \frac{b(s)^2}{\sigma(s)^2} (X^{(\alpha)}_s)^2
                   \, \dd s
             -\nu[X_t^{(\alpha)}]^2
         \right\}\right).
 \]
Heuristically, using \eqref{SEGED_Radon_Nicodym}, we have for all \ $\alpha$, $\beta\in\RR$,
 \ $\mu>0$, \ $\nu\geq 0$ \ and \ $t\in(0,T)$,
 \begin{align}\nonumber
   \Psi_t(\alpha,\mu,\nu)
       &=\EE\left[\exp\left\{-\mu\int_0^t\frac{b(s)^2}{\sigma(s)^2}(X^{(\beta)}_s)^2\,\dd s
                  -\nu[X_t^{(\beta)}]^2\right\}
                  \frac{\dd \PP_{X^{(\alpha)},\,t}}{\dd \PP_{X^{(\beta)},\,t}}
                    \left(X^{(\beta)}\big\vert_{[0,t]}\right)
                  \right] \\ \label{SEGED_DELYON_HU}
       &=\EE\left[\exp\left\{-\mu\int_0^t\frac{b(s)^2}{\sigma(s)^2}(X^{(\beta)}_s)^2\,\dd s
                  -\nu[X_t^{(\beta)}]^2
                  +(\alpha-\beta)\int_0^t\frac{b(s)}{\sigma(s)^2} X^{(\beta)}_s
                       \,\dd X_s^{(\beta)}\right.\right.\\ \nonumber
       &\phantom{=\EE\Big[\exp\Big\{\;}\left.\left.
                 -\frac{\alpha^2-\beta^2}{2}\int_0^t\frac{b(s)^2}{\sigma(s)^2}(X^{(\beta)}_s)^2\,\dd s
             \right\}\right].
 \end{align}
In what follows, using Theorem 1 in Delyon and Hu \cite{DelHu}, we give a precise derivation
 of \eqref{SEGED_DELYON_HU}.
For all \ $t\in(0,T),$ \ let \ $g:[0,t]\times\RR\to\RR,$ \ $h:[0,t]\times\RR\to\RR$ \ and
  \ $\sigma:[0,t]\times\RR\to\RR$ \ be defined by
 \begin{align*}
     g(u,x):=\alpha b(u)x,\qquad h(u,x):=(\beta-\alpha)\frac{b(u)}{\sigma(u)}x,
     \qquad \sigma(u,x):=\sigma(u),\quad \forall\;\,(u,x)\in[0,t]\times\RR.
 \end{align*}
Then \ $g$, \ $h$ \ and \ $\sigma$ \ are locally Lipschitz functions with respect
 to the second variable.
Let \ $f:C([0,t])\times C([0,t])\to\RR$,
 \[
  f(x,w):=\exp\left\{-\mu\int_0^t\frac{b(s)^2}{\sigma(s)^2}x(s)^2\,\dd s-\nu[x(t)]^2\right\},
  \qquad \forall\;\,(x,w)\in C([0,t])\times C([0,t]).
 \]
Using Theorem 1 in Delyon and Hu \cite{DelHu} with the above choices of \ $g$, \ $h$,
 \ $\sigma$ \ and \ $f$, \ we obtain for all \ $\alpha,$ $\beta\in\RR$, \ $\mu>0$,
 \ $\nu\geq 0$ \ and \ $t\in(0,T)$,
 \begin{align*}
   \Psi_t&(\alpha,\mu,\nu)
       =\EE\left[\exp\left\{-\mu\int_0^t\frac{b(u)^2}{\sigma(u)^2} (X^{(\beta)}_u)^2\,\dd u
                  -\nu[X_t^{(\beta)}]^2
                  -(\beta-\alpha)\int_0^t\frac{b(u)}{\sigma(u)} X^{(\beta)}_u
                       \,\dd B_u\right.\right.\\
       &\phantom{=\EE\Big[\exp\Big\{\qquad\qquad}\left.\left.
                 -\frac{(\beta-\alpha)^2}{2}\int_0^t\frac{b(u)^2}{\sigma(u)^2}(X^{(\beta)}_u)^2 \,\dd u
             \right\}\right].
 \end{align*}
By the SDE \eqref{special_SDE}, we conclude \eqref{SEGED_DELYON_HU}.

We check that for all \ $\beta\in\RR$ \ and \ $t\in[0,T)$, \
 \begin{align}\label{SEGED_LAPLACE1}
     \int_0^t\frac{b(s)}{\sigma(s)^2}X^{(\beta)}_s \,\dd X_s^{(\beta)}
        =\frac{1}{2}\left(\frac{b(t)}{\sigma(t)^2} (X^{(\beta)}_t)^2
          -\int_0^t\left[\frac{\dd}{\dd s}\left(\frac{b(s)}{\sigma(s)^2}\right)\right](X^{(\beta)}_s)^2\,\dd s
                 -\int_0^tb(s)\,\dd s\right).
 \end{align}
By It\^{o}'s rule (see, e.g., Liptser and Shiryaev \cite[Theorem 4.4]{LipShiI}), we get
 \begin{align}\label{SEGED_ITO}
  \begin{split}
  \dd\left(\frac{b(t)}{\sigma(t)^2}X^{(\beta)}_t\right)
     &=\left[\frac{\dd}{\dd t}\left(\frac{b(t)}{\sigma(t)^2}\right)\right]X^{(\beta)}_t\,\dd t
        +\frac{b(t)}{\sigma(t)^2}\,\dd X^{(\beta)}_t\\
     &=\left[\frac{\dd}{\dd t}\left(\frac{b(t)}{\sigma(t)^2}\right)\right]X^{(\beta)}_t\,\dd t
       +\beta\frac{b(t)^2}{\sigma(t)^2}X_t^{(\beta)}\,\dd t
       +\frac{b(t)}{\sigma(t)}\,\dd B_t,
       \quad t\in[0,T).
   \end{split}
 \end{align}
Now we verify that \ $(X_t^{(\beta)})_{t\in[0,T)}$ \ and
 \ $\left(\frac{b(t)}{\sigma(t)^2}X^{(\beta)}_t\right)_{t\in[0,T)}$
 \ are continuous semimartingales adapted to the filtration induced by \ $B$.
\ Consider the decomposition
 \[
  X_t^{(\beta)}
     =\exp\left\{\beta\int_0^t b(u)\,\dd u\right\}
      \int_0^t\sigma(s) \exp\left\{-\beta\int_0^s b(u)\,\dd u\right\}\,\dd B_s,
      \qquad t\in[0,T).
 \]
Here the deterministic function \ $\exp\left\{\beta\int_0^t b(u)\,\dd u\right\}$, \ $t\in[0,T)$,
 \ is monotone and hence has a finite variation over each finite interval of
 \ $[0,T)$, \ and then, by Jacod and Shiryaev \cite[Proposition 4.28, Chapter I]{JacShi},
 it is a semimartingale. Since
 \[
    \int_0^t\sigma(s) \exp\left\{-\beta\int_0^s b(u)\,\dd u\right\}\,\dd B_s,
      \qquad t\in[0,T),
 \]
 is a martingale with respect to the filtration induced by \ $B$,
 \ using Theorem 4.57 in Chapter I in Jacod and Shiryaev \cite{JacShi} with the function
 \ $f(x,y):=xy$, \ $x,y\in\RR$, \ we have \ $(X_t^{(\beta)})_{t\in[0,T)}$
 \ is a continuous semimartingale adapted to the filtration induced by \ $B$.
\ Similarly as above, using that by our assumptions, \ $\frac{b(t)}{\sigma(t)^2}$, \ $t\in[0,T)$,
 \ is continously differentiable, and hence has a finite variation over each finite interval of
 \ $[0,T)$, \ one can get
 \ $\left(\frac{b(t)}{\sigma(t)^2}X^{(\beta)}_t\right)_{t\in[0,T)}$
 \ is a continuous semimartingale adapted to the filtration induced by \ $B$.
\ Moreover, by \eqref{SEGED_ITO}, the cross-variation process of the continuous martingale
 parts of the processes \ $(X^{(\beta)}_t)_{t\in[0,T)}$ \ and
 \ $\left(\frac{b(t)}{\sigma(t)^2}X^{(\beta)}_t\right)_{t\in[0,T)}$ \ equals
 \[
    \int_0^t\sigma(s)\frac{b(s)}{\sigma(s)}\,\dd s
        =\int_0^t b(s)\,\dd s, \qquad t\in[0,T).
 \]
Hence, by integration by parts formula (see, e.g., Karatzas and Shreve \cite[page 155]{KarShr}),
 we have
 \begin{align*}
  \int_0^t\frac{b(s)}{\sigma(s)^2}X^{(\beta)}_s \,\dd X_s^{(\beta)}
     &=\frac{b(t)X^{(\beta)}_t}{\sigma(t)^2}X_t^{(\beta)}
       -\int_0^tX^{(\beta)}_s\,\dd \left(\frac{b(s)}{\sigma(s)^2}X^{(\beta)}_s\right)
       -\int_0^tb(s)\,\dd s\\
     &=\frac{b(t)}{\sigma(t)^2}(X^{(\beta)}_t)^2
       -\int_0^t\left[\frac{\dd}{\dd s}\left(\frac{b(s)}{\sigma(s)^2}\right)\right]
              (X^{(\beta)}_s)^2\,\dd s\\
     &\phantom{=\;}
       -\int_0^t\frac{b(s)}{\sigma(s)^2}X^{(\beta)}_s
              \,\dd X^{(\beta)}_s
       -\int_0^tb(s)\,\dd s,
       \qquad t\in[0,T),
 \end{align*}
 which gives us \eqref{SEGED_LAPLACE1}.

Then, using condition \eqref{Laplace_feltetel}, we have
 \begin{align}\label{SEGED_LAPLACE3}
  \begin{split}
     &\Psi_t(\alpha,\mu,\nu)\\
     &\;=\EE\left[\exp\left\{
                -\frac{1}{2}(2\mu+\alpha^2-\beta^2)\int_0^t\frac{b(s)^2}{\sigma(s)^2}
                           (X^{(\beta)}_s)^2\,\dd s
                -\frac{1}{2}\left(2\nu-\frac{(\alpha-\beta)b(t)}{\sigma(t)^2}\right)
                   (X_t^{(\beta)})^2\right.\right.\\
      &\phantom{\;=\EE\left[\exp\left\{\;\;\right.\right.}
            \left.\left.
       -\frac{\alpha-\beta}{2}\int_0^tb(s)\,\dd s
       -\frac{\alpha-\beta}{2}
          \int_0^t\left[\frac{\dd}{\dd s}\left(\frac{b(s)}{\sigma(s)^2}\right)\right]
           (X_s^{(\beta)})^2\,\dd s \right\}\right]\\
      &\;=\EE\left[\exp\left\{
              -\frac{1}{2}(2\mu+\alpha^2-\beta^2-2K(\alpha-\beta))\int_0^t\frac{b(s)^2}{\sigma(s)^2}
                            (X^{(\beta)}_s)^2\,\dd s\right.\right.\\
      &\phantom{\;=\EE\left[\exp\left\{\;\;\right.\right.}
            \left.\left.
                -\frac{1}{2}\left(2\nu-\frac{(\alpha-\beta)b(t)}{\sigma(t)^2}\right)(X^{(\beta)}_t)^2
               -\frac{\alpha-\beta}{2}\int_0^tb(s)\,\dd s
         \right\}\right].
   \end{split}
 \end{align}
We choose \ $\beta\in\RR$ \ such that \ $2\mu+\alpha^2-\beta^2-2K(\alpha-\beta)=0$.
\ Namely, let
 \[
   \beta:=K\pm\sqrt{2\mu+(\alpha-K)^2},
    \qquad \text{if}\quad \sign(b)=\pm\bone_{[0,T)}.
 \]
Then
 \begin{align}\label{SEGED_LAPLACE4}
  \Psi_t(\alpha,\mu,\nu)
     =\exp\left\{-\frac{\alpha-\beta}{2}\int_0^tb(s)\,\dd s\right\}
       \EE\left[\exp\left\{
                 -\frac{1}{2}\left(2\nu-\frac{(\alpha-\beta)b(t)}{\sigma(t)^2}\right)(X^{(\beta)}_t)^2
                 \right\}\right].
 \end{align}
The Laplace transform of a normally distributed random variable \ $\xi$ \ with mean \ $0$ \ and
 with variance \ $D>0$ \ is
 \begin{align}\label{SEGED3}
    \EE(\ee^{-s\xi^2})
       =\frac{1}{\sqrt{1+2sD}},
     \qquad s\geq 0.
 \end{align}
Since for all \ $t\in[0,T),$ \ $X_t^{(\beta)}$ \ is normally distributed with
 mean \ $0$ \ and with variance \ $V(t;\beta)$, \ using \eqref{SEGED3}
 we have for all \ $t\in[0,T)$,
 \[
    \EE\left[\exp\left\{-\frac{1}{2}\left(2\nu-(\alpha-\beta)\frac{b(t)}{\sigma(t)^2}\right)
                       (X^{(\beta)}_t)^2\right\}\right]
      =\frac{1}{\sqrt{1+\left(2\nu-(\alpha-\beta)\frac{b(t)}{\sigma(t)^2}\right)V(t;\beta)}}.
 \]
For this we have to check that
 \begin{align*}
   \frac{1}{2}\left(2\nu-(\alpha-\beta)\frac{b(t)}{\sigma(t)^2}\right)
       \geq 0, \qquad t\in[0,T).
 \end{align*}
This is satisfied, since \ $\nu\geq 0$ \ and for all \ $\alpha\in\RR$, \ $\mu>0$, \ we have
 \[
  \alpha-\beta
   =\alpha-K\mp\sqrt{2\mu+(\alpha-K)^2}
   =A_{\mu,\alpha,K}^{\pm},
 \]
 and hence \ $(\alpha-\beta)b(t)\leq 0$ \ for all \ $t\in[0,T)$ \ in both cases.
\proofend

\begin{Rem} \label{REM_technical}
Note that in Lemma \ref{Psi_general} we do not use the explicit solutions of the DE
 \eqref{Laplace_feltetel} given in Lemma \ref{LEMMA_DE}, since we wanted to demonstrate
  the role of condition \eqref{Laplace_feltetel} in the proof of Theorem \ref{Psi_general}.
By this condition, the process
 \ $\int_0^t \left[\frac{\dd}{\dd s}\left(\frac{b(s)}{\sigma(s)^2}\right)\right]
            (X^{(\beta)}_s)^2\,\dd s$,
 \ $t\in[0,T)$, \ has the form
 \ $-2K\int_0^t\frac{b(s)^2}{\sigma(s)^2}(X^{(\beta)}_s)^2\,\dd s$, \ $t\in[0,T)$,
 \ and hence $\int_0^t\frac{b(s)}{\sigma(s)^2}X^{(\beta)}_s\,\dd X^{(\beta)}_s$,
 \ can be expressed in terms of only the random variables
 \ $\big( X_t^{(\beta)} \big)^2$ \ and
 \ $\int_0^t \frac{b(s)^2}{\sigma(s)^2}(X_s^{(\beta)})^2 \, \dd s$,
 \ see formula \eqref{SEGED_LAPLACE1}.
As a consequence, in the calculation of \ $\Psi_t(\alpha,\mu,\nu)$ \ in the proof of
 Theorem \ref{Psi_general}, by the special choice of \ $\beta$, \ one can get
 rid of the stochastic integral
 \ $\int_0^t\frac{b(s)^2}{\sigma(s)^2}(X^{(\beta)}_s)^2\,\dd s$, \ see
 \eqref{SEGED_LAPLACE3} and \eqref{SEGED_LAPLACE4}.
\end{Rem}

In the next lemma we calculate explicitly the variance \ $V(t;\alpha)$ \ of
 \ $X_t^{(\alpha)}$ \ for all \ $t\in[0,T)$.

\begin{Lem}\label{LEMMA_VARIANCE2}
Let \ $(X_t^{(\alpha)})_{t\in[0,T)}$ \ be the process given by the SDE \eqref{special_SDE},
 where \ $b$ \ is given by \eqref{Laplace_feltetel2}.
Then
 \begin{align*}
   V(t;\alpha)
    =\begin{cases}
       \frac{C}{\alpha-K}\left(B_{K,C}(t)^\alpha-B_{K,C}(t)^K\right)
           & \quad \text{if \ $\alpha\ne K$,}\\[2mm]
       C B_{K,C}(t)^K \ln(B_{K,C}(t))
           & \quad \text{if \ $\alpha=K$,}
     \end{cases}
 \end{align*}
 where \ $B_{K,C}(t)$, \ $t\in[0,T)$, \ is defined in Theorem \ref{THEOREM_LAPLACE}.
\end{Lem}

\noindent{\bf Proof.}
First let us suppose that \ $b(t)>0$ \ for all \ $t\in[0,T)$.
\ Then \ $C$ \ is positive, since by \ $b(0)>0$, \ $K\int_0^0\sigma(u)^2\,\dd u+C$
 \ should be positive.
If \ $\alpha\ne K$ \ and \ $K\ne0$, \ by \eqref{SEGED17_uj}, we have
 for all \ $t\in[0,T)$, \
 \begin{align*}
   V(t;\alpha)
     &=\int_0^t
         \left(\frac{K\int_0^t\sigma(u)^2\,\dd u+C}
            {K\int_0^s\sigma(u)^2\,\dd u+C}\right)^{\frac{\alpha}{K}}\sigma(s)^2\,\dd s\\
     &=\frac{1}{K-\alpha}\left(K\int_0^t\sigma(u)^2\,\dd u+C\right)^{\frac{\alpha}{K}}
         \left(\left(K\int_0^t\sigma(u)^2\,\dd u+C\right)^{\frac{K-\alpha}{K}}
                \!-C^{\frac{K-\alpha}{K}}\right),
 \end{align*}
 which yields the assertion in case of \ $\alpha\ne K$, \ $K\ne0$.

The other cases can be handled similarly.

Let us suppose now that \ $b(t)<0$ \ for all \ $t\in[0,T)$. \
For all \ $\beta\in\RR$, \ let us consider the process \ $(N_t^{(\beta)})_{t\in[0,T)}$ \
 given by the SDE
 \begin{align*}
  \left\{
   \begin{aligned}
    \dd N_t^{(\beta)}
    &= \beta \, \widetilde b(t) N_t^{(\beta)} \, \dd t + \sigma(t) \, \dd B_t ,
    \qquad t \in [0,T) , \\
    N_0^{(\beta)} &= 0 ,
   \end{aligned}
  \right.
 \end{align*}
 where \ $\widetilde b(t):=-b(t)$, \ $t\in[0,T)$. \
Then, by uniqueness of a strong solution, the process \ $(X_t^{(\alpha)})_{t\in[0,T)}$
 \ given by the SDE \eqref{special_SDE} and the process \ $(N_t^{(-\alpha)})_{t\in[0,T)}$
 \ coincide and hence \ $V(t;\alpha)=V_{N^{(-\alpha)}}(t)$, \ $t\in[0,T)$,
 \ where \ $V_{N^{(-\alpha)}}(t):=\EE(N^{(-\alpha)}_t)^2$, \ $t\in[0,T)$.
\ Moreover, \ $V_{N^{(-\alpha)}}(t)$, \ $t\in[0,T)$, \ is given by the formulae in the present
 Lemma \ref{LEMMA_VARIANCE2} with \ $(\alpha,K,C)$ \ is replaced by \ $(-\alpha,-K,-C)$.
\ Since these formulae are invariant under the above defined replacement, we have the assertion.
\proofend

\noindent{\bf Proof of Theorem \ref{THEOREM_LAPLACE}.}
First we check that for all \ $K\in\RR$,
 \begin{align}\label{SEGED_THEOREM_LAPLACE1}
  \int_0^t b(s)\,\dd s = \frac{1}{2}\ln(B_{K,C}(t)),
     \qquad t\in[0,T).
 \end{align}
If \ $K\ne0$, \ then
 \begin{align*}
   \int_0^t b(s)\,\dd s
    &=\int_0^t \frac{\sigma(s)^2}{2\left(K\int_0^s\sigma(u)^2\,\dd u+C\right)}\,\dd s
     =\frac{1}{2K}\ln\left(K\int_0^t\sigma(u)^2\,\dd u+C\right)
      -\frac{1}{2K}\ln C\\
    &=\frac{1}{2}\ln\left(1+\frac{K}{C}\int_0^t\sigma(u)^2\,\dd u\right)^{\frac{1}{K}}
     =\frac{1}{2}\ln(B_{K,C}(t)), \qquad t\in[0,T),
 \end{align*}
 and if \ $K=0$, \ then \ $\int_0^t b(s)\,\dd s=\int_0^t \frac{\sigma(s)^2}{2C}\,\dd s=\frac{1}{2}\ln(B_{K,C}(t))$,
 \ $t\in[0,T)$.
\ By Lemmas \ref{Psi_general} and \ref{LEMMA_VARIANCE2}, using also \eqref{SEGED_THEOREM_LAPLACE1},
 for all \ $\mu>0$, \ $\nu\geq 0$, \ and \ $t\in[0,T)$, \ we have
 \begin{align*}
  \psi_t(\alpha,\mu,\nu)
  &=\left(\frac{\exp\left\{-\frac{A_{\mu,\alpha,K}^{\pm}}{2}\ln(B_{K,C}(t))\right\}}
        {1+\left(2\nu-\frac{A_{\mu,\alpha,K}^{\pm}}{2\left(K\int_0^t\sigma(s)^2\,\dd s+C\right)}\right)
         \frac{C}{\pm\sqrt{2\mu+(\alpha-K)^2}}
          \left(B_{K,C}(t)^{K\pm\sqrt{2\mu+(\alpha-K)^2}}-B_{K,C}(t)^K\right)}
      \right)^{\frac{1}{2}}\\[2mm]
  &=\left(\frac{B_{K,C}(t)^{-\frac{\alpha-K}{2}}}
            {D}
       \right)^{\frac{1}{2}},
  \end{align*}
 where
 \begin{align*}
  D&:= B_{K,C}(t)^{\mp\frac{\sqrt{2\mu+(\alpha-K)^2}}{2}}\\
    &\phantom{:=\,}
     +\frac{4\nu\left(K\int_0^t\sigma(s)^2\,\dd s+C\right)-A_{\mu,\alpha,K}^{\pm}}
                  {\pm2\sqrt{2\mu+(\alpha-K)^2}}
               \left( B_{K,C}(t)^{\pm\frac{\sqrt{2\mu+(\alpha-K)^2}}{2}}
                      -B_{K,C}(t)^{\mp\frac{\sqrt{2\mu+(\alpha-K)^2}}{2}} \right)\\
    &=\left(\frac{1}{2}
                 \pm\frac{4\nu\left(K\int_0^t\sigma(s)^2\,\dd s+C\right)-\alpha+K}
                         {2\sqrt{2\mu+(\alpha-K)^2}}\right)
                     B_{K,C}(t)^{\pm\frac{\sqrt{2\mu+(\alpha-K)^2}}{2}}\\
    &\phantom{:=\,}
      +\left(\frac{1}{2}
                 \mp\frac{4\nu\left(K\int_0^t\sigma(s)^2\,\dd s+C\right)-\alpha+K}
                         {2\sqrt{2\mu+(\alpha-K)^2}}\right)
                     B_{K,C}(t)^{\mp\frac{\sqrt{2\mu+(\alpha-K)^2}}{2}},
 \end{align*}
 which yields the assertion.
\proofend

\begin{Rem}
Note that formula \eqref{LAPLACE} in Lemma \ref{Psi_general} for the joint Laplace transform
 of \eqref{SEGED_abstract1} depends on the sign of the function \ $\sign(b)$,
 \ but in Theorem \ref{THEOREM_LAPLACE} it turned out that the sign is indifferent.
We also remark that the case \ $b(t)<0$, \ $t\in[0,T)$, \ can be traced back to the case
 \ $b(t)>0$, \ $t\in[0,T)$, \ using the same arguments that are written for the case
 \ $b(t)<0$, \ $t\in[0,T)$, \ at the end of the proof of Lemma \ref{LEMMA_VARIANCE2}.
The point is that the formulae in Theorem \ref{THEOREM_LAPLACE} are invariant
 under the replacement of \ $(\alpha,b,K,C)$ \ with \ $(-\alpha,-b,-K,-C)$.
\end{Rem}

In the next two remarks we consider special cases of Theorem \ref{THEOREM_LAPLACE}.

\begin{Rem}
As a special case of Theorem \ref{THEOREM_LAPLACE}, one can get back formula \eqref{LP_OU_FORMULA}
 due to Liptser and Shiryaev \cite[Lemma 17.3]{LipShiII}, and also the well-known
 Cameron--Martin formula for a standard Wiener process.
Namely, let \ $T:=\infty$, \ $b(t):=1$, \ $t\geq 0$, \ and \ $\sigma(t):=1$,
 \ $t\geq 0$.
\ Let us consider the process \ $(X^{(\alpha)}_t)_{t\in[0,T)}$ \ given by the SDE
 \eqref{special_SDE}, which is the usual Ornstein--Uhlenbeck process starting from \ $0$.
\ Clearly, \ $\frac{\dd}{\dd t}\left(\frac{b(t)}{\sigma(t)^2}\right)=0$, \ $t>0$,
 and hence Theorem \ref{THEOREM_LAPLACE} with \ $\nu=0$, \ $K=0$ \ and with
 \ $C=\frac{1}{2}$ \ implies \eqref{LP_OU_FORMULA}.
With \ $\alpha=0$, \ we get back the Cameron--Martin formula for a standard Wiener process,
  \[
     \EE\exp\left\{-\mu\int_0^t(B_u)^2\,\dd u\right\}
          =\frac{1}{\sqrt{\cosh(t\sqrt{2\mu})}},
          \qquad t\geq 0,\quad \mu>0,
 \]
 see, e.g., Liptser and Shiryaev \cite[formula (7.147)]{LipShiI}.
\end{Rem}

\begin{Rem}\label{REM_MANSUY}
Let \ $T\in(0,\infty)$, \ $b(t):=-\frac{1}{T-t}$, \ $t\in[0,T)$, \ and \ $\sigma(t):=1$,
 \ $t\in[0,T)$. \
Let us consider the process \ $(X^{(\alpha)}_t)_{t\in[0,T)}$ \ given by the SDE
 \eqref{special_SDE}.
Hence condition \eqref{Laplace_feltetel2} is satisfied with \ $K:=\frac{1}{2}$ \ and
 \ $C:=-\frac{T}{2}$, \ and clearly, \ $B_{K,C}(t)=(1-t/T)^2$,
 \ $t\in[0,T$).
\ Then Theorem \ref{THEOREM_LAPLACE} with \ $\nu=0$ \ and \ $\alpha=0$ \ implies that for all
 \ $\mu>0$ \ and \ $t\in[0,T)$,
  \begin{align*}
      &\EE\exp\left\{-\frac{\mu}{2}\int_0^t\frac{(B_u)^2}{(T-u)^2}\,\dd u\right\}\\
      &\qquad\quad
       =\frac{\left(1-\frac{t}{T}\right)^{\frac{1}{4}}}
         {\sqrt{\cosh\left(\ln\left(1-\frac{t}{T}\right)\sqrt{\mu+\frac{1}{4}}\right)
           +\frac{1}{2\sqrt{\mu+\frac{1}{4}}}
             \sinh\left(\ln\left(1-\frac{t}{T}\right)\sqrt{\mu+\frac{1}{4}}\right)}}.
 \end{align*}
An easy calculation shows that for all \ $\mu>0$ \ and \ $t\in[0,T)$,
 \begin{align*}
   \EE\exp\left\{-\frac{\mu}{2}\int_0^t\frac{(B_u)^2}{(T-u)^2}\,\dd u\right\}
    =\frac{\left(\frac{T-t}{T}\right)^{\frac{1+\sqrt{4\mu+1}}{4}}}
     {\sqrt{1-\frac{1+\sqrt{4\mu+1}}{2\sqrt{4\mu+1}}\left(1-\left(1-\frac{t}{T}\right)^{\sqrt{4\mu+1}}\right)}}.
 \end{align*}
This is the corrected formula of Proposition 5 in Mansuy \cite{Man}, which contains a
 misprint.
\end{Rem}

\section{Maximum likelihood estimation via Laplace transform}\label{Section_Maximum_Likelihood}

As a special case of \eqref{SEGED_Radon_Nicodym}, the measures \ $\PP_{X^{(\alpha)},\,t}$
 \ and \ $\PP_{X^{(0)},\,t}$ \ are equivalent for all \ $\alpha\in\RR$ \ and for all \ $t\in(0,T)$,
 \ and
 \begin{align*} 
     \frac{\dd \PP_{X^{(\alpha)},\,t}}{\dd \PP_{X^{(0)},\,t}}
       \left(X^{(\alpha)}\big\vert_{[0,t]}\right)
        = \exp \left\{ \alpha
                  \int_0^t
                   \frac{b(s)}{\sigma(s)^2}X_s^{(\alpha)} \, \dd X_s^{(\alpha)}
                  -\frac{\alpha^2}{2}
                   \int_0^t
                    \frac{b(s)^2}{\sigma(s)^2}(X_s^{(\alpha)})^2 \, \dd s \right\}.
 \end{align*}
Here \ $\PP_{X^{(0)},\,t}$ \ is nothing else but the Wiener measure on
 \ $\big(C([0,t]),\cB(C([0,t]))\big)$.

For all \ $t\in(0,T)$, \ the maximum likelihood estimator
 \ $\halpha_t$ \ of the parameter \ $\alpha$ \ based on the
 observation \ $(X_s^{(\alpha)})_{s \in [0,\,t]}$ \ is defined by
 \[
   \halpha_t
    :=\argmax_{\alpha\in\RR}
       \ln\left(\frac{\dd \PP_{X^{(\alpha)},\,t}}{\dd \PP_{X^{(0)},\,t}}
                 \left(X^{(\alpha)}\big\vert_{[0,t]}\right)\right).
 \]

The following lemma due to Barczy and Pap \cite[Lemma 1]{BarPap} guarantees the
 existence of a unique MLE of \ $\alpha$.

\begin{Lem}\label{LEMMA_SPEC_FELTETEL}
For all \ $\alpha\in\RR$ \ and \ $t\in(0,T)$, \ we have
 \[
   \PP\left( \int_0^t
               \frac{b(s)^2}{\sigma(s)^2}(X_s^{(\alpha)})^2 \, \dd s
               > 0 \right)
   =1.
 \]
\end{Lem}

By Lemma \ref{LEMMA_SPEC_FELTETEL}, for all \ $t\in(0,T),$ \ there exists a
 unique maximum likelihood estimator \ $\halpha_t$ \ of the
 parameter \ $\alpha$ \ based on the observation \ $(X_s^{(\alpha)})_{s \in [0,\,t]}$
 \ given by
 \[
   \halpha_t
   = \frac{\int_0^t
            \frac{b(s)}{\sigma(s)^2}X_s^{(\alpha)} \, \dd X_s^{(\alpha)}}
          {\int_0^t
            \frac{b(s)^2}{\sigma(s)^2}(X_s^{(\alpha)})^2 \, \dd s},
   \qquad  t\in(0,T).
 \]
To be more precise, by Lemma \ref{LEMMA_SPEC_FELTETEL}, for all
 \ $t\in(0,T),$ \ the MLE \ $\halpha_t$ \ exists $\PP$-almost surely.
Using the SDE \eqref{special_SDE} we obtain
 \begin{align}\label{SEGED20_uj}
   \halpha_t - \alpha
   = \frac{\int_0^t
            \frac{b(s)}{\sigma(s)}X_s^{(\alpha)} \, \dd B_s}
          {\int_0^t
            \frac{b(s)^2}{\sigma(s)^2}(X_s^{(\alpha)})^2 \, \dd s},
      \qquad t\in(0,T).
 \end{align}
For all \ $t\in(0,T)$, \ the Fisher information for \ $\alpha$ \ contained in
 the observation \ $(X_s^{(\alpha)})_{s\in[0,\,t]}$, \ is defined by
 \[
   I_{\alpha}(t)
   := \EE \left( \frac{\partial}{\partial\alpha}
                 \ln \left(\frac{\dd \PP_{X^{(\alpha)},\,t}}{\dd \PP_{X^{(0)},\,t}}
                  \left( X^{(\alpha)} \big \vert_{[0,\,t]} \right) \right) \right)^2
   = \int_0^t
      \frac{b(s)^2}{\sigma(s)^2}\,\EE \big( X_s^{(\alpha)} \big)^2  \, \dd s ,
 \]
 where the last equality follows by the SDE \eqref{special_SDE} and Karatzas and Shreve
 \cite[Proposition 3.2.10]{KarShr}.
Note that, by the conditions on \ $b$ \ and \ $\sigma$,
 \ $I_{\alpha} : (0,T) \to (0,\infty)$ \ is an increasing function.
\ Now we calculate the Fisher information \ $I_{\alpha}(t)$, \ $t\in(0,T)$,
 \ explicitly.

\begin{Lem}\label{LEMMA_FISHER2}
Let \ $(X_t^{(\alpha)})_{t\in[0,T)}$ \ be the process given by the SDE \eqref{special_SDE},
 where \ $b$ \ is given by \eqref{Laplace_feltetel2}.
Then for all \ $t\in(0,T)$,
 \begin{align*}
   I_{\alpha}(t)
     =\begin{cases}
       \frac{1}{4(\alpha-K)^2}(B_{K,C}(t)^{\alpha-K}-1)
        -\frac{1}{4(\alpha-K)}\ln(B_{K,C}(t))
           & \quad \text{if \ $\alpha\ne K$,} \\[2mm]
       \frac{1}{8}(\ln(B_{K,C}(t)))^2
                      & \quad \text{if \ $\alpha=K$,}
      \end{cases}
 \end{align*}
 where \ $B_{K,C}(t)$, \ $t\in[0,T)$, \ is defined in Theorem \ref{THEOREM_LAPLACE}.
\end{Lem}

\noindent{\bf Proof.}
First let us suppose that \ $b(t)>0$ \ for all \ $t\in[0,T)$.
\ Then \ $C$ \ is positive, since by \ $b(0)>0$, \ $K\int_0^0\sigma(u)^2\,\dd u+C$
 \ should be positive.
In case of \ $\alpha\ne K$ \ and \ $K\ne0$, \ by Lemma \ref{LEMMA_VARIANCE2},
 we get for all \ $t\in(0,T)$, \
 \begin{align*}
   I_{\alpha}(t)
      &=\int_0^t\frac{\sigma(s)^2}{4\left(K\int_0^s\sigma(u)^2\,\dd u+C\right)^2}
                V(s;\alpha)\,\dd s
       =\int_0^t \frac{\sigma(s)^2}{4C(\alpha-K)}\big(B_{K,C}(s)^{\alpha-2K}-B_{K,C}(s)^{-K}\big)\,\dd s \\
      &=\int_0^t\frac{\sigma(s)^2}{4C(\alpha-K)}
        \left(\left(1+\frac{K}{C}\int_0^s\sigma(u)^2\,\dd u\right)^{\frac{\alpha-2K}{K}}
                     -\left(1+\frac{K}{C}\int_0^s\sigma(u)^2\,\dd u\right)^{-1}
               \right)\,\dd s,
 \end{align*}
 which yields the assertion in case of \ $\alpha\ne K$ \ and \ $K\ne0$.

The other cases can be handled similarly.

The case \ $b(t)<0$, \ $t\in[0,T)$, \ can be handled similarly to what is written
 for the case \ $b(t)<0$, \ $t\in[0,T)$, \ at the end of the proof of Lemma \ref{LEMMA_VARIANCE2}.
The point is that the formulae in the present Lemma \ref{LEMMA_FISHER2} are invariant under
 the replacement of \ $(\alpha,b,K,C)$ \ with \ $(-\alpha,-b,-K,-C)$.
\hspace*{5cm} \proofend

Later on we intend to prove limit theorems for the MLE \ $\widehat\alpha_t$ \ of \ $\alpha$
 \ normalized by Fisher information \ $I_{\alpha}(t)$.
\ For proving these limit theorems, condition \ $\lim_{t\uparrow T}I_{\alpha}(t)=\infty$
 \ plays a crucial role.
In what follows we examine under what additional conditions on \ $b$ \ and \ $\sigma$,
 \ $\lim_{t\uparrow T}I_{\alpha}(t)=\infty$ \ is satisfied.

\begin{Lem}\label{LEMMA_FISHER_ASYMPTOTIC2}
Let \ $(X_t^{(\alpha)})_{t\in[0,T)}$ \ be the process given by the SDE \eqref{special_SDE},
 where \ $b$ \ is given by \eqref{Laplace_feltetel2}.
In case of \ $K\ne 0$, \
 \[
   \lim_{t\uparrow T}I_{\alpha}(t)=\infty
     \qquad \Longleftrightarrow \qquad
     \lim_{t\uparrow T}\int_0^t\sigma(u)^2\,\dd u
       =\begin{cases}
         \infty & \quad \text{if \ $\frac{C}{K}>0$,}\\
         -\frac{C}{K} & \quad \text{if \ $\frac{C}{K}<0$.}
         \end{cases}
 \]
In case of \ $K=0$, \ we have \ $\lim_{t\uparrow T}I_{\alpha}(t)=\infty$ \ holds
 if and only if \ $\lim_{t\uparrow T}\int_0^t\sigma(u)^2\,\dd u=\infty$.
\end{Lem}

\noindent{\bf Proof.}
First we note that \ $C\ne0$, \ since by \ $b(0)\ne0$, \ $K\int_0^0\sigma(u)^2\,\dd u+C$ \
 should be not zero.
Now we check that for all \ $K\in\RR$,
 \begin{align}\label{SEGED_FISHER}
    \lim_{t\uparrow T}I_{\alpha}(t)=\infty
     \qquad \Longleftrightarrow \qquad
    \lim_{t\uparrow T} B_{K,C}(t)^{\alpha-K}\in\{0,\infty\}.
 \end{align}
If \ $\alpha\ne K$, \ by Lemma \ref{LEMMA_FISHER2}, we get
 \begin{align*}
   I_\alpha(t)
     &=\frac{1}{4(\alpha-K)^2}
      \Big(\exp\big\{(\alpha-K)\ln(B_{K,C}(t))\big\}-(\alpha-K)\ln(B_{K,C}(t))-1\Big)\\
     &=\frac{1}{4(\alpha-K)^2}f\left(\ln(B_{K,C}(t)^{\alpha-K})\right),
 \end{align*}
 where \ $f(x):=\ee^x-x-1,$ \ $x\in\RR$.
\ Using that the function \ $\int_0^t\sigma(u)^2\,\dd u$, \ $t\in[0,T)$, \ is monotone increasing,
 we have \ $\lim_{t\uparrow T}B_{K,C}(t)$ \ exists. Hence
 \begin{align*}
    \lim_{t\uparrow T}I_{\alpha}(t)=\infty
     \qquad \Longleftrightarrow \qquad
    \lim_{t\uparrow T} \ln (B_{K,C}(t)^{\alpha-K})\in\{-\infty,\infty\},
 \end{align*}
 which implies \eqref{SEGED_FISHER}.
A similar argument shows that \eqref{SEGED_FISHER} is valid also in case of \ $\alpha=K$.
\ Hence, by the definition of \ $B_{K,C}(t)$, \ we have in case of \ $K\ne0$,
 \[
    \lim_{t\uparrow T}I_{\alpha}(t)=\infty
     \qquad \Longleftrightarrow \qquad
    \lim_{t\uparrow T}\left(1+\frac{K}{C}\int_0^t\sigma(s)^2\,\dd s\right)^{\frac{\alpha-K}{K}}
      \in\{0,\infty\},
 \]
 and in case of \ $K=0$, \
 \[
    \lim_{t\uparrow T}I_{\alpha}(t)=\infty
     \qquad \Longleftrightarrow \qquad
    \lim_{t\uparrow T}\exp\left\{\frac{\alpha}{C}\int_0^t\sigma(s)^2\,\dd s\right\}
      \in\{0,\infty\}.
 \]
This implies the assertion.
\proofend

Note that if the function \ $b:[0,T)\to\RR\setminus\{0\}$ \ is given by \eqref{Laplace_feltetel2}
 and if we suppose also that \ $K\ne0$, \ $\frac{C}{K}<0$, \ then,
 by Lemma \ref{LEMMA_FISHER_ASYMPTOTIC2}, we have
 \begin{align}\label{SEGED_FISHER3}
    C=-K\lim_{t\uparrow T}\int_0^t\sigma(u)^2\,\dd u
     =:-K\int_0^T\sigma(u)^2\,\dd u\in\RR\setminus\{0\},
 \end{align}
 and hence
 \begin{align*}
    b(t)=\frac{\sigma(t)^2}{2\left(K\int_0^t\sigma(u)^2\,\dd u-K\int_0^T\sigma(u)^2\,\dd u\right)}
        =\frac{\sigma(t)^2}{-2K\int_t^T\sigma(u)^2\,\dd u},
       \qquad t\in[0,T),
 \end{align*}
 which is nothing else but the form \eqref{SEGED_abstract2} of \ $b$.
\ Moreover, by Lemma \ref{LEMMA_FISHER_ASYMPTOTIC2},
 we have \ $\lim_{t\uparrow T}I_{\alpha}(t)=\infty$ \ holds in this case.

In all what follows we will suppose that the function \ $b$ \ is given
 by \eqref{SEGED_abstract2} with some \ $K\ne0$, \ where \ $\int_0^T\sigma(u)^2\,\dd u<\infty$,
 \ and in this case, as an application of the explicit form of the joint Laplace transform
 of \eqref{SEGED_abstract1}, we will give a complete description of the asymptotic
 behavior of the MLE \ $\widehat\alpha_t$ \ of \ $\alpha$ \ as \ $t\uparrow T$.
\ In the other cases (for which \ $\lim_{t\uparrow T}I_{\alpha}(t)=\infty$) \
 the asymptotic behavior of the MLE \ $\widehat\alpha_t$ \ as \ $t\uparrow T$
 \ may be worked out using the same arguments as follows, but we do not consider these cases.

For our later purposes, we examine the asymptotic behavior of
 \ $I_{\alpha}(t)$ \ as \ $t\uparrow T$.

\begin{Lem}\label{LEMMA_FISHER3}
Let \ $(X_t^{(\alpha)})_{t\in[0,T)}$ \ be the process given by the SDE \eqref{special_SDE},
 where \ $b$ \ is given by \eqref{SEGED_abstract2}
 with some \ $K\ne0$ \ and we suppose that \ $\int_0^T\sigma(s)^2\,\dd s<\infty$.
Then in case of \ $\sign(\alpha-K)=-\sign(K)$,
 \[
   \lim_{t\uparrow T}\frac{I_{\alpha}(t)}
      {\frac{1}{4(K-\alpha)^2}
        \left(\frac{\int_0^T\sigma(s)^2\,\dd s}
              {\int_t^T\sigma(s)^2\,\dd s}\right)^{\frac{K-\alpha}{K}}}
       =1,
 \]
 in case of \ $\alpha=K$,
 \[
   \lim_{t\uparrow T}\frac{I_{\alpha}(t)}
      {\frac{1}{8K^2}
        \left(\ln\left(\int_t^T\sigma(s)^2\,\dd s\right)\right)^2}
       =1,
 \]
 and in case of \ $\sign(\alpha-K)=\sign(K)$,
 \[
   \lim_{t\uparrow T}\frac{I_{\alpha}(t)}
      {\frac{1}{4K(K-\alpha)}
         \ln\left(\int_t^T\sigma(s)^2\,\dd s\right)}
       =1.
 \]
\end{Lem}

The next lemma is about the asymptotic behavior of the Laplace transform of the denominator
 in \eqref{SEGED20_uj}.

\begin{Lem}\label{LEMMA_FISHER4}
Let \ $(X_t^{(\alpha)})_{t\in[0,T)}$ \ be the process given by the SDE \eqref{special_SDE},
 where \ $b$ \ is given by \eqref{SEGED_abstract2}
 with some \ $K\ne0$ \ and we suppose that \ $\int_0^T\sigma(s)^2\,\dd s<\infty$.
Then
 \begin{align}\label{SEGED48}
   \frac{1}{I_{\alpha}(t)}
   \int_0^t \frac{b(u)^2}{\sigma(u)^2}(X_u^{(\alpha)})^2 \, \dd u
   \distr
   \begin{cases}
    (W_1)^2 & \text{if \ $\sign(\alpha-K)=-\sign(K)$,} \\[2mm]
    2\int_0^1 (W_s)^2 \, \dd s & \text{if \ $\alpha = K$,} \\[2mm]
    1 & \text{if \ $\sign(\alpha-K)=\sign(K)$,}
   \end{cases}
 \end{align}
 as \ $t \uparrow T$, \ where \ $(W_s)_{s\in[0,1]}$ \ is a standard Wiener process.
In fact, in case of \ $\alpha = K$, \ for all \ $t \in (0,T)$,
 \begin{align}\label{SEGED51}
   \frac{1}{I_{K}(t)}
   \int_0^t \frac{b(u)^2}{\sigma(u)^2} (X_u^{(K)})^2 \, \dd u
   \distre
   2\int_0^1 (W_s)^2 \, \dd s ,
   \qquad t\in(0,T).
 \end{align}
\end{Lem}

\noindent{\bf Proof.}
We will show that for all \ $\mu>0$,
 \begin{align}\label{SEGED_DENOM_LAPLACE1}
    \lim_{t \uparrow T}
    \EE \exp \left\{ - \frac{\mu}{I_{\alpha}(t)}
                       \int_0^t \frac{b(u)^2}{\sigma(u)^2}(X_u^{(\alpha)})^2 \, \dd u \right\}
    = \begin{cases}
       \frac{1}{\sqrt{1 + 2\mu}} & \text{if \ $\sign(\alpha-K)=-\sign(K)$,} \\[2mm]
       \frac{1}{\sqrt{\cosh(2\sqrt{\mu})}} & \text{if \ $\alpha = K$,}
        \\[2mm]
       \ee^{-\mu} & \text{if \ $\sign(\alpha-K)=\sign(K)$.}
      \end{cases}
 \end{align}
In fact, in case of \ $\alpha = K$, \ we prove that for all \ $t \in (0,T)$ \ and \ $\mu \geq 0$,
 \begin{align}\label{SEGED_DENOM_LAPLACE2}
   \EE \exp \left\{ - \frac{\mu}{I_{K}(t)}
                      \int_0^t \frac{b(u)^2}{\sigma(u)^2} (X_u^{(K)})^2 \, \dd u \right\}
   = \frac{1}{\sqrt{\cosh(2\sqrt{\mu})}}.
 \end{align}

First we suppose that \ $K<0$. \ Then we have \ $b(t)>0$, \ $t\in[0,T)$,
 \ and the function \ $b$ \ satisfies the DE \eqref{Laplace_feltetel}.
By \eqref{SEGED_FISHER3},
 \begin{align}\label{SEGED52}
   B_{K,C}(t)
      =\left(\frac{\int_t^T\sigma(s)^2\,\dd s}{\int_0^T\sigma(s)^2\,\dd s}\right)^{\frac{1}{K}},
      \qquad t\in[0,T),\quad K\ne0,
 \end{align}
 and hence, by Theorem \ref{THEOREM_LAPLACE}, for all \ $\alpha\in\RR$, \ $\mu>0$ \ and \ $t\in(0,T)$,
 \ we get
 \begin{align}\label{SEGED_LAPLACE2}
   \EE \exp \left\{ - \frac{\mu}{I_{\alpha}(t)}
                      \int_0^t \frac{b(u)^2}{\sigma(u)^2}(X_u^{(\alpha)})^2 \, \dd u \right\}
        =\frac{1}{\sqrt{C_{\mu,\alpha,K}(t)}}
 \end{align}
 where
 \begin{align*}
   C_{\mu,\alpha,K}(t)
   &:=\left(\frac{1}{2}+\frac{\alpha-K}
                 {2\widetilde A_{\mu,\alpha,K}(t)}\right)
           \left(\frac{\int_t^T\sigma(s)^2\,\dd s}
                  {\int_0^T\sigma(s)^2\,\dd s}\right)^{\frac{\alpha-K-\widetilde A_{\mu,\alpha,K}(t)}{2K}} \\
     &\phantom{:=\;\;}
       +\left(\frac{1}{2}-\frac{\alpha-K}
                 {2\widetilde A_{\mu,\alpha,K}(t)}\right)
           \left(\frac{\int_t^T\sigma(s)^2\,\dd s}
                  {\int_0^T\sigma(s)^2\,\dd s}\right)^{\frac{\alpha-K+\widetilde A_{\mu,\alpha,K}(t)}{2K}},
 \end{align*}
 and
 \begin{align*}
   \widetilde A_{\mu,\alpha,K}(t):=\sqrt{\frac{2\mu}{I_{\alpha}(t)}+(\alpha-K)^2}.
 \end{align*}

Now we consider the case \ $K<0$ \ and \ $\alpha>K$.
\ Using that \ $\lim_{t \uparrow T}I_{\alpha}(t)=\infty$ \ and \ $\alpha-K>0$,
 \ we have \ $\lim_{t \uparrow T}\widetilde A_{\mu,\alpha,K}(t)=\alpha-K$.
\ Then, using Lemma \ref{LEMMA_FISHER3} and that \ $\lim_{x\downarrow 0}x^x=1$,
 \ an easy calculation shows that
 \begin{align*}
   \lim_{t \uparrow T}
     \left(\frac{1}{2}+\frac{\alpha-K}{2\widetilde A_{\mu,\alpha,K}(t)}\right)
     &  \left(\frac{\int_t^T\sigma(s)^2\,\dd s}
                  {\int_0^T\sigma(s)^2\,\dd s}\right)^{\frac{\alpha-K-\widetilde A_{\mu,\alpha,K}(t)}{2K}}\\
     &=\lim_{t \uparrow T}
       \left(\left(\frac{\int_t^T\sigma(s)^2\,\dd s}
          {\int_0^T\sigma(s)^2\,\dd s}\right)
        ^{-1+\sqrt{8\mu\left(\frac{\int_0^T\sigma(s)^2\,\dd s}
          {\int_t^T\sigma(s)^2\,\dd s}\right)^{\frac{\alpha-K}{K}}+1}}
          \right)^{\frac{-\alpha+K}{2K}}\\
     &=\lim_{t \uparrow T}
              \left(\frac{\int_t^T\sigma(s)^2\,\dd s}
                    {\int_0^T\sigma(s)^2\,\dd s}\right)
               ^{\frac{8\mu\left(\frac{\int_t^T\sigma(s)^2\,\dd s}
                          {\int_0^T\sigma(s)^2\,\dd s}\right)^{\frac{-\alpha+K}{K}}}
                  {1+\sqrt{8\mu\left(\frac{\int_t^T\sigma(s)^2\,\dd s}
                          {\int_0^T\sigma(s)^2\,\dd s}\right)^{\frac{-\alpha+K}{K}}+1}}
                  \frac{-\alpha+K}{2K}}
     =1.
 \end{align*}
\ Moreover,
 \begin{align*}
    &\lim_{t \uparrow T}
       \left(\frac{1}{2}-\frac{\alpha-K}{2\widetilde A_{\mu,\alpha,K}(t)}\right)
       \left(\frac{\int_t^T\sigma(s)^2\,\dd s}
                  {\int_0^T\sigma(s)^2\,\dd s}\right)^{\frac{\alpha-K+\widetilde A_{\mu,\alpha,K}(t)}{2K}}\\[2mm]
      &=\lim_{t \uparrow T}
            \left(\frac{1}{2}-\frac{1}{2\sqrt{8\mu\left(\frac{\int_0^T\sigma(s)^2\,\dd s}
                          {\int_t^T\sigma(s)^2\,\dd s}\right)^{\frac{\alpha-K}{K}}+1}}\right)
            \!\!\left(\frac{\int_t^T\sigma(s)^2\,\dd s}
                  {\int_0^T\sigma(s)^2\,\dd s}\right)
            ^{\frac{\alpha-K}{2K}
               \left(1+\sqrt{8\mu\left(\frac{\int_0^T\sigma(s)^2\,\dd s}
                          {\int_t^T\sigma(s)^2\,\dd s}\right)^{\frac{\alpha-K}{K}}+1}\right)}\\[2mm]
    &=\lim_{t \uparrow T}
      \frac{4\mu\left(\frac{\int_t^T\sigma(s)^2\,\dd s}
                  {\int_0^T\sigma(s)^2\,\dd s}\right)
            ^{\frac{\alpha-K}{2K}
               \left(-1+\sqrt{8\mu\left(\frac{\int_t^T\sigma(s)^2\,\dd s}
                          {\int_0^T\sigma(s)^2\,\dd s}\right)^{\frac{-\alpha+K}{K}}+1}\right)}}
        {\sqrt{8\mu\left(\frac{\int_t^T\sigma(s)^2\,\dd s}
                          {\int_0^T\sigma(s)^2\,\dd s}\right)^{\frac{-\alpha+K}{K}}+1}
        \left(1+\sqrt{8\mu\left(\frac{\int_t^T\sigma(s)^2\,\dd s}
                          {\int_0^T\sigma(s)^2\,\dd s}\right)^{\frac{-\alpha+K}{K}}+1}\right)}
    =2\mu,
 \end{align*}
 since the denominator tends to \ $2$ \ as \ $t\uparrow T$, \ and
 \begin{align*}
   \lim_{t \uparrow T}
      \left(\frac{\int_t^T\sigma(s)^2\,\dd s}
                  {\int_0^T\sigma(s)^2\,\dd s}\right)
            &^{\frac{\alpha-K}{2K}
               \left(-1+\sqrt{8\mu\left(\frac{\int_t^T\sigma(s)^2\,\dd s}
                          {\int_0^T\sigma(s)^2\,\dd s}\right)^{\frac{-\alpha+K}{K}}+1}\right)}\\
   &=\lim_{t \uparrow T}
      \left(\frac{\int_t^T\sigma(s)^2\,\dd s}
                  {\int_0^T\sigma(s)^2\,\dd s}\right)
        ^{\frac{\alpha-K}{2K}
          \frac{8\mu\left(\frac{\int_t^T\sigma(s)^2\,\dd s}
                  {\int_0^T\sigma(s)^2\,\dd s}\right)^{\frac{-\alpha+K}{K}}}
          {1+\sqrt{8\mu\left(\frac{\int_t^T\sigma(s)^2\,\dd s}
                          {\int_0^T\sigma(s)^2\,\dd s}\right)^{\frac{-\alpha+K}{K}}+1}}}
   =1.
 \end{align*}
Hence, by \eqref{SEGED_LAPLACE2}, we have \eqref{SEGED_DENOM_LAPLACE1} in case of \ $K<0$
 \ and \ $\alpha>K$.
\ By \eqref{SEGED3}, for all \ $\mu>0$, \ we have
 \[
  \EE(\ee^{-\mu(W_1)^2})=\frac{1}{\sqrt{1+2\mu}},
 \]
 and the unicity of Laplace transform implies \eqref{SEGED48} in case
 of \ $K<0$ \ and \ $\alpha>K$.

Now we consider the case \ $K<0$ \ and \ $\alpha=K$.
\ For all \ $t\in(0,T)$ \ and \ $\mu>0$, \ by \eqref{SEGED_LAPLACE2}, we get
 \begin{align*}
   \EE\exp\left\{
                -\frac{\mu}{I_{K}(t)}\int_0^t\frac{b(s)^2}{\sigma(s)^2} (X_s^{(K)})^2\,\dd s
           \right\}
     &=\frac{1}
       {\sqrt{\frac{1}{2}
           \left(\frac{\int_t^T\sigma(s)^2\,\dd s}{\int_0^T\sigma(s)^2\,\dd s}\right)
                ^{\sqrt{\frac{\mu}{2K^2 I_{K}(t)}}}
        +\frac{1}{2}
           \left(\frac{\int_t^T\sigma(s)^2\,\dd s}{\int_0^T\sigma(s)^2\,\dd s}\right)
             ^{-\sqrt{\frac{\mu}{2K^2I_{K}(t)}}}}}\\
     &=\frac{1}
       {\sqrt{\frac{1}{2}\ee^{-2\sqrt{\mu}}
         +\frac{1}{2}\ee^{2\sqrt{\mu}}}}
      =\frac{1}{\sqrt{\cosh(2\sqrt{\mu})}},
 \end{align*}
 where the last but one equality follows from the fact, by Lemma \ref{LEMMA_FISHER2},
 in case of \ $K<0$ \ and \ $\alpha=K$ \ we have
 \begin{align}\label{SEGED_FISHER2}
  \sqrt{I_{K}(t)}
      =\frac{1}{2\sqrt{2}K}\ln\left(\frac{\int_t^T\sigma(s)^2\,\dd s}{\int_0^T\sigma(s)^2\,\dd s}\right),
     \qquad t\in(0,T),
 \end{align}
 and from the fact that \ $x^{\frac{1}{\ln x}}=\ee$ \ for all \ $x>0$.
\ By formula (1.9.3) in Borodin and Salminen \cite[Part II, Section 1]{BorSal}, we get
 \begin{align*}
   \EE\exp\left\{-2\mu\int_0^1(W_u)^2\,\dd u\right\}
      =\frac{1}{\sqrt{\cosh(2\sqrt{\mu})}},
    \qquad \mu>0,
 \end{align*}
 and the unicity of Laplace transform implies \eqref{SEGED_DENOM_LAPLACE2} and \eqref{SEGED51}
 in case of \ $K<0$ \ and \ $\alpha=K$.

Now we consider the case \ $K<0$ \ and \ $\alpha<K$.
\ Using that \ $\lim_{t\uparrow T}I_{\alpha}(t)=\infty$, \
 we have \ $\lim_{t\uparrow T}\widetilde A_{\mu,\alpha,K}(t)=-(\alpha-K)$, \ since \ $\alpha-K<0$.
\ Then
 \begin{align*}
    \lim_{t \uparrow T}
      \left(\frac{1}{2}+\frac{\alpha-K}{2\widetilde A_{\mu,\alpha,K}(t)}\right)
        \left(\frac{\int_t^T\sigma(s)^2\,\dd s}
                   {\int_0^T\sigma(s)^2\,\dd s}\right)^{\frac{\alpha-K-\widetilde A_{\mu,\alpha,K}(t)}{2K}}
        =0.
 \end{align*}
Moreover, by Lemma \ref{LEMMA_FISHER3}, we get
 \begin{align*}
    \lim_{t \uparrow T}
       &\left(\frac{1}{2}-\frac{\alpha-K}{2\widetilde A_{\mu,\alpha,K}(t)}\right)
        \left(\frac{\int_t^T\sigma(s)^2\,\dd s}
                  {\int_0^T\sigma(s)^2\,\dd s}\right)^{\frac{\alpha-K+\widetilde A_{\mu,\alpha,K}(t)}{2K}}\\
    &\;=\lim_{t \uparrow T}
            \left(\frac{\int_t^T\sigma(s)^2\,\dd s}
                   {\int_0^T\sigma(s)^2\,\dd s}\right)^
          {\frac{\alpha-K+\sqrt{\frac{8\mu}{\frac{1}{K(K-\alpha)}
                       \ln\left(\int_t^T\sigma(s)^2\,\dd s\right)}+(\alpha-K)^2}}{2K}}\\
    &\;=\lim_{t \uparrow T}
          \left(\frac{\int_t^T\sigma(s)^2\,\dd s}
                   {\int_0^T\sigma(s)^2\,\dd s}\right)^
          {-\frac{K-\alpha}{2K}
           \left(1-\sqrt{\frac{8\mu K}{(K-\alpha)\ln\left(\int_t^T\sigma(s)^2\,\dd s\right)}+1}\right)}\\
    &\;=\lim_{t \uparrow T}
                         \left(\frac{\int_t^T\sigma(s)^2\,\dd s}
                   {\int_0^T\sigma(s)^2\,\dd s}\right)^
                 {\frac{K-\alpha}{2K}
                   \frac{\frac{8\mu K}{(K-\alpha)\ln\left(\int_t^T\sigma(s)^2\,\dd s\right)}}
                   {1+\sqrt{\frac{8\mu K}{(K-\alpha)\ln\left(\int_t^T\sigma(s)^2\,\dd s\right)}+1}}}
    =\ee^{2\mu},
 \end{align*}
 since
  \[
     \lim_{t \uparrow T}
        \left(\frac{\int_t^T\sigma(s)^2\,\dd s}
                   {\int_0^T\sigma(s)^2\,\dd s}\right)^
            {\frac{1}{\ln\left(\int_t^T\sigma(s)^2\,\dd s\right)}}
            =\ee.
  \]
Hence, by \eqref{SEGED_LAPLACE2} and the unicity of Laplace transform, we have \eqref{SEGED_DENOM_LAPLACE1}
 and \eqref{SEGED48} in case of \ $K<0$ \ and \ $\alpha<K$.

The case \ $K>0$ \ can be handled in the same way as at the end of the proof of Lemma \ref{LEMMA_FISHER2}.
\proofend

\begin{Thm}\label{THM_alpha_1/2_general}
Let \ $(X_t^{(\alpha)})_{t\in[0,T)}$ \ be the process given by the SDE \eqref{special_SDE},
 where \ $b$ \ is given by \eqref{SEGED_abstract2}
 with some \ $K\ne0$ \ and we suppose that \ $\int_0^T\sigma(s)^2\,\dd s<\infty$.
Then
 \[
   \sqrt{I_{\alpha}(t)} \left(\widehat\alpha_t-\alpha\right)
    \distr
   \begin{cases}
     \cN(0,1)  & \text{if \ $\sign(\alpha-K)=\sign(K)$,}\\[1mm]
     -\frac{\sign(K)}{\sqrt{2}} \, \frac{\int_0^1 W_s \, \dd W_s}{\int_0^1 (W_s)^2 \, \dd s}
              & \text{if \ $\alpha=K$,}
    \end{cases}
 \]
 as \ $t\uparrow T$, \ where \ $(W_s)_{s\in[0,1]}$ \ is a standard Wiener process.
In fact, in case of \ $\alpha=K$, \ for all \ $t\in(0,T)$,
 \begin{align}\label{SEGED_SINGULAR1}
    \sqrt{I_{K}(t)}\left(\widehat\alpha_t-K\right)
       \distre
        -\frac{\sign(K)}{2\sqrt{2}}\frac{(W_1)^2-1}{\int_0^1(W_s)^2\,\dd s}
        =
         -\frac{\sign(K)}{\sqrt{2}}\frac{\int_0^1W_s\,\dd W_s}{\int_0^1(W_s)^2\,\dd s}.
 \end{align}
\end{Thm}

\noindent{\bf Proof.}
First we suppose that \ $K<0$.
\ Then we have \ $b(t)>0$, \ $t\in[0,T)$, \ and the function \ $b$ \ satisfies the DE
 \eqref{Laplace_feltetel}.
By the SDE \eqref{special_SDE} and \eqref{SEGED_LAPLACE1}, we have for all
 \ $\alpha\in\RR$ \ and \ $t\in[0,T)$,
 \begin{align}\label{SEGED_max_lik_uj3}
  \begin{split}
   \int_0^t\frac{b(s)}{\sigma(s)}X^{(\alpha)}_s\,\dd B_s
     &= \int_0^t\frac{b(s)}{\sigma(s)^2}X^{(\alpha)}_s\,\dd X^{(\alpha)}_s
       -\alpha \int_0^t\frac{b(s)^2}{\sigma(s)^2}(X^{(\alpha)}_s)^2\,\dd s\\
     &= \frac{b(t)}{2\sigma(t)^2}(X^{(\alpha)}_t)^2
       -\frac{1}{2}\int_0^tb(s)\,\dd s
       -\left(\alpha-K\right)
        \int_0^t\frac{b(s)^2}{\sigma(s)^2}(X^{(\alpha)}_s)^2\,\dd s.
  \end{split}
 \end{align}

Now let us suppose that \ $K<0$ \ and \ $\alpha<K$.
\ By Lemma \ref{LEMMA_FISHER3}, \ $\lim_{t\uparrow T}I_{\alpha}(t)=\infty$ \ holds,
 and Lemma \ref{LEMMA_FISHER4} implies that
 \begin{align*}
     \frac{1}{I_{\alpha}(t)}\int_0^t\frac{b(s)^2}{\sigma(s)^2}(X^{(\alpha)}_s)^2\,\dd s
         \stoch 1
      \qquad \text{as \ $t\uparrow T$,}
 \end{align*}
 where \ $\stoch$ \ denotes convergence in probability.
Indeed, if \ $K<0$ \ and\ $\alpha<K$, \ then the limit in \eqref{SEGED48} is \ $1$,
 \ which is a constant, and hence convergence in distribution implies convergence in probability.
Hence we can apply Theorem 4 in Barczy and Pap \cite{BarPap} with
 \ $Q(t):=\frac{1}{\sqrt{I_{\alpha}(t)}}$, \ $t\in(0,T)$,
 \ and \ $\eta:=1$, \ and then we have the assertion in case of \ $K<0$ \ and
 \ $\alpha<K$.

Now let us suppose that \ $K<0$ \ and \ $\alpha=K$.
\ By \eqref{SEGED20_uj} and \eqref{SEGED_max_lik_uj3}, we get
 \begin{align*}
   \widehat\alpha_t - K
     =&\frac{\frac{b(t)}{2\sigma(t)^2}\big(X^{(K)}_t\big)^2
           -\frac{1}{2}\int_0^tb(s)\,\dd s}
         {\int_0^t\frac{b(s)^2}{\sigma(s)^2}(X_s^{(K)})^2\,\dd s},
    \qquad t\in(0,T).
 \end{align*}
Then for all \ $t\in(0,T)$,
 \begin{align*}
   \sqrt{I_{K}(t)}&\left(\widehat\alpha_t-K\right)
      =\frac{1}{2\sqrt{2}}
        \frac{\frac{1}{\sqrt{2I_{K}(t)}}\frac{b(t)}{\sigma(t)^2}\big(X^{(K)}_t\big)^2
           -\frac{1}{\sqrt{2I_{K}(t)}}\int_0^tb(s)\,\dd s}
         {\frac{1}{2I_{K}(t)}\int_0^t\frac{b(s)^2}{\sigma(s)^2}(X_s^{(K)})^2\,\dd s}.
 \end{align*}
To prove \eqref{SEGED_SINGULAR1}, it is enough to check that
 \begin{align}\label{SEGED_SINGULAR2}
   \begin{split}
    &\left(\frac{1}{\sqrt{2I_{K}(t)}}\frac{b(t)}{\sigma(t)^2}\big(X^{(K)}_t\big)^2,
          \frac{1}{2I_{K}(t)}\int_0^t\frac{b(s)^2}{\sigma(s)^2}(X_s^{(K)})^2\,\dd s\right)\\
     &\quad\distre\left((W_1)^2,\int_0^1(W_s)^2\,\dd s\right)
      \distre\left(1+2\int_0^1 W_s\,\dd W_s,\int_0^1(W_s)^2\,\dd s\right),
       \qquad t\in(0,T),
   \end{split}
 \end{align}
 and
 \begin{align}\label{SEGED_SINGULAR3}
      \int_0^tb(s)\,\dd s=\sqrt{2I_{K}(t)},
        \qquad t\in(0,T).
 \end{align}
Using that for all \ $\mu>0$ \ and \ $\nu\geq 0$,
 \[
  \EE\exp\left\{-\mu\int_0^1(W_s)^2\,\dd s-\nu[W_1]^2\right\}
        =\frac{1}{\sqrt{\cosh(\sqrt{2\mu})+\frac{2\nu}{\sqrt{2\mu}}\sinh(\sqrt{2\mu})}},
 \]
 (see, e.g., formula (1.9.3) in Borodin and Salminen \cite[Part II, Section 1]{BorSal},
 or as a special case of our Theorem \ref{THEOREM_LAPLACE}), to prove the first equality in distribution
 of \eqref{SEGED_SINGULAR2}, it is enough to verify that for all \ $\mu>0$ \ and \ $\nu\geq 0$,
 \begin{align*}
   \EE\exp&\left\{
                -\frac{\mu}{2I_{K}(t)}\int_0^t\frac{b(s)^2}{\sigma(s)^2}(X_s^{(K)})^2\,\dd s
                -\frac{\nu}{\sqrt{2I_{K}(t)}}\frac{b(t)}{\sigma(t)^2}(X^{(K)}_t)^2
           \right\}\\
          &=\frac{1}{\sqrt{\cosh(\sqrt{2\mu})+\frac{2\nu}{\sqrt{2\mu}}\sinh(\sqrt{2\mu})}},
           \qquad t\in(0,T).
 \end{align*}
By Theorem \ref{THEOREM_LAPLACE}, we get for all \ $t\in(0,T)$,
 \begin{align*}
    \EE&\exp\left\{
         -\frac{\mu}{2I_{K}(t)}\int_0^t\frac{b(s)^2}{\sigma(s)^2}(X_s^{(K)})^2\,\dd s
         -\frac{\nu}{\sqrt{2I_{K}(t)}}\frac{b(t)}{\sigma(t)^2}(X^{(K)}_t)^2
            \right\}\\
     &=\frac{1}
       {\sqrt{\left(\frac{1}{2}-\frac{\nu}{\sqrt{2\mu}}\right)
           \left(\frac{\int_t^T\sigma(s)^2\,\dd s}{\int_0^T\sigma(s)^2\,\dd s}\right)
                 ^{\sqrt{\frac{\mu}{4K^2I_{K}(t)}}}
        +\left(\frac{1}{2}+\frac{\nu}{\sqrt{2\mu}}\right)
           \left(\frac{\int_t^T\sigma(s)^2\,\dd s}{\int_0^T\sigma(s)^2\,\dd s}\right)
                 ^{-\sqrt{\frac{\mu}{4K^2I_{K}(t)}}}}}\\
     &=\frac{1}
       {\sqrt{\left(\frac{1}{2}-\frac{\nu}{\sqrt{2\mu}}\right)\ee^{-\sqrt{2\mu}}
         +\left(\frac{1}{2}+\frac{\nu}{\sqrt{2\mu}}\right)\ee^{\sqrt{2\mu}}}}
      =\frac{1}{\sqrt{\cosh(\sqrt{2\mu})+\frac{2\nu}{\sqrt{2\mu}}\sinh(\sqrt{2\mu})}},
 \end{align*}
 where the last but one equality follows from \eqref{SEGED_FISHER2}
 and from the fact that \ $x^{\frac{1}{\ln x}}=\ee$ \ for all \ $x>0$.
\ Hence, by the uniqueness of Laplace transform, for all \ $t\in(0,T)$,
 \ the joint distribution of
 \[
   \frac{1}{2I_{K}(t)}
     \int_0^t \frac{b(s)^2}{\sigma(s)^2}(X_s^{(K)})^2\,\dd s
   \qquad\text{and}\qquad
   \frac{1}{\sqrt{2I_{K}(t)}}\frac{b(t)}{\sigma(t)^2}(X_t^{(K)})^2
 \]
 is the same as the joint distribution of \ $\int_0^1(W_s)^2\,\dd s$ \ and
 \ $(W_1)^2$.
\ Finally, by It\^o's formula,
 \[
    \int_0^1 W_s \, \dd W_s=\frac{1}{2}((W_1)^2-1),
 \]
 and hence for all \ $t \in (0,T)$, \ we have \eqref{SEGED_SINGULAR2}.
We note that \ $\frac{\int_0^1 W_s \, \dd W_s}{\int_0^1 (W_s)^2 \, \dd s}$ \ is the
 limit distribution of the Dickey--Fuller statistic, see, e.g., the Ph.D. thesis of Bobkoski
 \cite{Bob}, or (7.14) and Theorem 9.5.1 in Tanaka \cite{Tan}.

\noindent Now we check \eqref{SEGED_SINGULAR3}.
Since \ $K<0$ \ and \ $\alpha=K$, \ using \eqref{SEGED_FISHER2}, we get for all
 \ $t\in(0,T)$,
 \begin{align*}
     \int_0^tb(s)\,\dd s
      =\int_0^t\frac{\sigma(s)^2}{-2K\int_s^T\sigma(u)^2\,\dd u}\,\dd s
      =\frac{1}{2K}
         \ln\left(\frac{\int_t^T\sigma(u)^2\,\dd u}{\int_0^T\sigma(u)^2\,\dd u}\right)
      =\sqrt{2I_{K}(t)}.
 \end{align*}

Let us suppose now that \ $K>0$.
\ Then \ $b(t)<0$ \ for all \ $t\in[0,T)$.
\ The statement in this case can be obtained from the case \ $b(t)>0$ \ for all \ $t\in[0,T)$,
 \ using the arguments at the end of the proof of Lemma \ref{LEMMA_VARIANCE2}.
The point is that we need to consider the replacement of \ $(\alpha,b,K)$ \ with
 \ $(-\alpha,-b,-K)$ \ and, with the notations introduced in the proof of Lemma \ref{LEMMA_VARIANCE2},
 to take into account that
 \ $\widehat{(-\alpha)}_t^{(N^{(-\alpha)})}=-\widehat\alpha_t^{(X^{(\alpha)})}$, \ $t\in(0,T)$.
\proofend

\begin{Rem}
We note that Theorem \ref{THM_alpha_1/2_general} can be derived from our more
 general results, namely, from Barczy and Pap \cite[Theorems 5 and 10]{BarPap}.
We also remark that using these results one can also weaken the conditions on
 \ $b$ \ and \ $\sigma$ \ in Theorem \ref{THM_alpha_1/2_general}.
\end{Rem}

\begin{Rem}\label{REM_CAUCHY}
In case of \ $\sign(\alpha-K)=-\sign(K)$, \ under the conditions of Theorem \ref{THM_alpha_1/2_general},
 one can prove that
 \[
   \sqrt{I_{\alpha}(t)} \left(\widehat\alpha_t-\alpha\right)
    \distr\zeta
    \text{\qquad as \ $t\uparrow T,$}
 \]
 where \ $\zeta$ \ is a standard Cauchy distributed random variable,
 see, e.g., Luschgy \cite[Section 4.2]{Lus1} or Barczy and Pap \cite{BarPap}.
The proof in this case is based on a martingale limit theorem, and we do not
 know whether one can find a proof using the explicit form of the joint Laplace transform
 of \eqref{SEGED_abstract1}. Lemma \ref{LEMMA_FISHER4} implies only
 \begin{align}\label{SEGED50}
    \frac{1}{I_{\alpha}(t)}
   \int_0^t \frac{b(u)^2}{\sigma(u)^2}(X_u^{(\alpha)})^2 \, \dd u
   \distr \cN(0,1)^2
   \qquad \text{as} \quad t \uparrow T.
 \end{align}
However, using a martingale limit theorem, one can prove that the convergence
 in \eqref{SEGED50} holds almost surely (with some appropriate random variable
 \ $\xi^2$ \ as the limit).
To be able to use Theorem 4 in Barczy and Pap \cite{BarPap},
 we need convergence in probability in \eqref{SEGED50}.
Hence the question is whether we can improve the convergence in distribution in \eqref{SEGED50}
 to convergence in probability using only the explicit form of the joint Laplace transform
 of \eqref{SEGED_abstract1}.
We do not know if one can find such a technique.
\end{Rem}

The next theorem is about the (asymptotic) behavior of the MLE \ of
 \ $\alpha=K$, \ $K\ne0$ \ using an appropriate {\sl random} normalizing factor.

\begin{Thm}\label{THM_alpha_1/2_general_random}
Let \ $(X_t^{(K)})_{t\in[0,T)}$ \ be the process given by the SDE \eqref{special_SDE},
 where \ $b$ \ is given by \eqref{SEGED_abstract2}
 with some \ $K\ne0$ \ and we suppose that \ $\int_0^T\sigma(s)^2\,\dd s<\infty$.
\ Then for all \ $t\in(0,T)$,
 \begin{align*}
   \left(\!\int_0^t\frac{b(u)^2}{\sigma(u)^2} (X_u^{(K)})^2 \,\dd u\!\right)^{\frac{1}{2}}
     \!\left(\widehat\alpha_t-K\right)
     \distre -\sign(K)
      \frac{\int_0^1 W_u\,\dd W_u}{\left(\int_0^1 (W_u)^2\,\dd u\right)^{\frac{1}{2}}}
     =-\frac{\sign(K)}{2}\frac{(W_1)^2-1}{\left(\int_0^1 (W_u)^2\,\dd u\right)^{\frac{1}{2}}}.
 \end{align*}
\end{Thm}

\noindent{\bf Proof.}
First we suppose that \ $K<0$.
\ By \eqref{SEGED_SINGULAR1} and \eqref{SEGED_SINGULAR2}, we have for all \ $\alpha\in\RR$
 \ and for all \ $t\in(0,T)$,
 \begin{align*}
    &\left(\int_0^t\frac{b(u)^2}{\sigma(u)^2}(X_u^{(K)})^2\,\dd u\right)^{\frac{1}{2}}
        \left(\widehat\alpha_t-K\right)
     =\sqrt{I_{K}(t)}(\widehat\alpha_t-K)
       \left(\frac{1}{I_{K}(t)}
       \int_0^t\frac{b(u)^2}{\sigma(u)^2}(X_u^{(K)})^2 \,\dd u\right)^{\frac{1}{2}}\\
    &\qquad\qquad\qquad
     \distre\frac{1}{\sqrt{2}}\frac{\int_0^1W_u\,\dd W_u}{\int_0^1(W_u)^2\,\dd u}
              \left(2\int_0^1 (W_u)^2 \, \dd u\right)^{\frac{1}{2}}
      =\frac{\int_0^1W_u\,\dd W_u}{\left(\int_0^1 (W_u)^2 \, \dd u\right)^{\frac{1}{2}}},
       \qquad t\in(0,T),
 \end{align*}
 which implies the assertion using It\^o's formula.

The case $K>0$ can be handled in the same way as at the end of the proof of Theorem
 \ref{THM_alpha_1/2_general}. \proofend

\begin{Rem}\label{REM_BARCZY_PAP}
We note that, by Barczy and Pap \cite[Corollaries 9 and 11]{BarPap},
 under the conditions \ $\int_0^T\sigma(s)^2\,\dd s<\infty$ \ and \eqref{SEGED_abstract2},
 we have for all \ $\alpha\ne K$, \ $K\ne0$, \ the MLE of \ $\alpha$ \ is asymptotically normal
 with a corresponding {\sl random} normalizing factor, namely, for all \ $\alpha\ne K$, \ $K\ne0$,
 \begin{align*}
    \left(\int_0^t\frac{b(u)^2}{\sigma(u)^2}(X_u^{(\alpha)})^2 \,\dd u\right)^{\frac{1}{2}}
        \left(\widehat\alpha_t-\alpha\right)
        \distr\cN(0,1)
        \qquad \text{as\quad $t\uparrow T$.}
 \end{align*}
\end{Rem}

As a consequence of Theorem \ref{THM_alpha_1/2_general_random}, giving an illuminating counterexample,
 we show that Remark 1.47 in Prakasa Rao \cite{Rao} contains a mistake.

\begin{Rem}\label{REM_RAO}
By giving a counterexample, we show that
 condition (1.5.26) in Remark 1.47 in Prakasa Rao \cite{Rao} is not
 enough to assure (1.5.35) in Prakasa Rao \cite{Rao}.
By \eqref{SEGED20_uj}, we have for all \ $\alpha\in\RR$ \ and \ $t\in(0,T),$ \
 \begin{align}\label{SEGED_RAO_MISTAKE}
    \left(\int_0^t\frac{b(u)^2}{\sigma(u)^2}(X_u^{(\alpha)})^2 \,\dd u\right)^{\frac{1}{2}}
        (\widehat\alpha_t-\alpha)
       =\frac{\frac{1}{\sqrt{I_{\alpha}(t)}}
         \int_0^t\frac{b(u)}{\sigma(u)}X_u^{(\alpha)}\,\dd B_u}
         {\left(\frac{1}{I_{\alpha}(t)}
         \int_0^t\frac{b(u)^2}{\sigma(u)^2}(X_u^{(\alpha)})^2 \,\dd u\right)^{1/2}}.
 \end{align}
By Lemma \ref{LEMMA_FISHER4} (under its conditions), we have
 \[
   \frac{1}{I_{K}(t)}
   \int_0^t \frac{b(u)^2}{\sigma(u)^2}(X_u^{(K)})^2 \, \dd u
   \distre
   2\int_0^1 (W_u)^2 \, \dd u ,
   \qquad t\in(0,T).
 \]
Hence if Remark 1.47 in Prakasa Rao \cite{Rao} were true, then we would have
 \begin{align*}
      &\left(\frac{1}{\sqrt{I_{K}(t)}}\int_0^t\frac{b(s)}{\sigma(s)}X_s^{(K)}\,\dd B_s,
            \frac{1}{I_{K}(t)}\int_0^t\frac{b(s)^2}{\sigma(s)^2}(X_s^{(K)})^2\,\dd s\right)\\
     &\qquad\qquad
       \distr\left(\left(2\int_0^1 (W_u)^2\,\dd u\right)^{\frac{1}{2}}\xi\,,\;2\int_0^1 (W_u)^2\,\dd u\right)
        \quad \text{as \ $t\uparrow T,$}
 \end{align*}
 where \ $\xi$ \ is a standard normally distributed random variable independent of
 \ $\int_0^1 (W_u)^2\,\dd u.$
\ By \eqref{SEGED_RAO_MISTAKE} and continuous mapping theorem, we would have
 \begin{align*}
   \left(\int_0^t\frac{b(u)^2}{\sigma(u)^2}(X_u^{(K)})^2\,\dd u\right)^{\frac{1}{2}}
        (\widehat\alpha_t-K)
     \distr \frac{\left(2\int_0^1 (W_u)^2\,\dd u\right)^{\frac{1}{2}}\xi}
                 {\left(2\int_0^1 (W_u)^2\,\dd u\right)^{\frac{1}{2}}}
             =\xi \qquad \text{as \ $t\uparrow T$,}
 \end{align*}
 which is a contradiction, since, by Theorem \ref{THM_alpha_1/2_general_random},
 the limit distribution is
 \[
    -\frac{\sign(K)}{2}\frac{(W_1)^2-1}{\left(\int_0^1 (W_u)^2 \, \dd u\right)^{\frac{1}{2}}}.
 \]
Note that this limit distribution can not be a standard normal distribution, see, e.g.,
 Feigin \cite[Section 2]{Fei2}. Indeed, in case of \ $K<0$,
 \begin{align*}
   \PP\left(-\frac{\sign(K)}{2}
     \frac{(W_1)^2-1}{\left(\int_0^1 (W_u)^2 \, \dd u\right)^{\frac{1}{2}}}>0\right)
     &=\PP((W_1)^2>1)=2(1-\Phi(1)),
 \end{align*}
 which is not equal to \ $\PP(\cN(0,1)>0)=\frac{1}{2}$.
\ In case of \ $K>0$, \ we can arrive at a contradiction similarly.
\end{Rem}

The next theorem is about the strong consistency of the MLE of \ $\alpha$. \

\begin{Thm}\label{THM_CONSISTENCY}
Let \ $(X_t^{(\alpha)})_{t\in[0,T)}$ \ be the process given by the SDE \eqref{special_SDE},
 where \ $b$ \ is given by \eqref{SEGED_abstract2}
 with some \ $K\ne0$ \ and we suppose that \ $\int_0^T\sigma(s)^2\,\dd s<\infty$.
 \ Then the maximum likelihood estimator of \ $\alpha$ \ is strongly consistent,
 i.e., for all \ $\alpha\in\RR$,
 \[
   \PP\Big(\lim_{t\uparrow T}\widehat\alpha_t=\alpha\Big)=1.
 \]
\end{Thm}

\noindent{\bf Proof.}
First we suppose that \ $K<0$.
\ Then we have \ $b(t)>0$, $t\in[0,T)$, \ and the function \ $b$ \ satisfies the DE
 \eqref{Laplace_feltetel}.
We check that for all \ $\alpha\in\RR$,
 \[
   \EE \exp \left\{-\lim_{t\uparrow T}\int_0^t \frac{b(u)^2}{\sigma(u)^2}(X_u^{(\alpha)})^2 \,\dd u
             \right\}
   =\lim_{t\uparrow T}
      \EE\exp\left\{-\int_0^t\frac{b(u)^2}{\sigma(u)^2} (X_u^{(\alpha)})^2\,\dd u
               \right\}
   =0.
 \]
The first equality follows from monotone convergence theorem, and the second one can
 be derived as follows.
Using \eqref{SEGED52} and Theorem \ref{THEOREM_LAPLACE} with \ $\mu:=1$ \ and \ $\nu:=0$,
 \ we get for all \ $t\in(0,T)$,
 \begin{align*}
   \EE\exp\left\{-\int_0^t \frac{b(u)^2}{\sigma(u)^2}(X_u^{(\alpha)})^2 \, \dd u \right\}
        =\frac{1}{\sqrt{C_{\alpha,K}(t)}}
 \end{align*}
 where
 \begin{align*}
   C_{\alpha,K}(t)
    &:=\left(\frac{1}{2}+\frac{\alpha-K}
                 {2\sqrt{2+(\alpha-K)^2}}\right)
           \left(\frac{\int_t^T\sigma(s)^2\,\dd s}
                  {\int_0^T\sigma(s)^2\,\dd s}\right)^{\frac{-K+\alpha-\sqrt{2+(\alpha-K)^2}}{2K}}\\
     &\phantom{:=\;\;}
       +\left(\frac{1}{2}-\frac{\alpha-K}{2\sqrt{2+(\alpha-K)^2}}\right)
           \left(\frac{\int_t^T\sigma(s)^2\,\dd s}
                  {\int_0^T\sigma(s)^2\,\dd s}\right)^{\frac{-K+\alpha+\sqrt{2+(\alpha-K)^2}}{2K}}.
 \end{align*}

\noindent In case of \ $\alpha-K\geq 0$, \ we have \ $\sqrt{2+(\alpha-K)^2}>\alpha-K$ \ and hence
  \begin{align}\label{SEGED_CONSISTENCY1}
     &\lim_{t\uparrow T}
         \left(\frac{\int_t^T\sigma(s)^2\,\dd s}
                   {\int_0^T\sigma(s)^2\,\dd s}\right)^
            {\frac{-K+\alpha-\sqrt{2+(\alpha-K)^2}}{2K}}
         =0,\\\label{SEGED_CONSISTENCY2}
     &\lim_{t\uparrow T}
         \left(\frac{\int_t^T\sigma(s)^2\,\dd s}
                   {\int_0^T\sigma(s)^2\,\dd s}\right)^
            {\frac{-K+\alpha+\sqrt{2+(\alpha-K)^2}}{2K}}
         =\infty.
  \end{align}

\noindent In case of \ $\alpha-K<0$, \ we have \ $\sqrt{2+(\alpha-K)^2}>-(\alpha-K)$
 \ and hence \eqref{SEGED_CONSISTENCY1} and \eqref{SEGED_CONSISTENCY2} are satisfied again.
 Since
 \[
   \frac{1}{2}-\frac{\alpha-K}{2\sqrt{2+(\alpha-K)^2}}
     =\frac{\sqrt{2+(\alpha-K)^2}-\alpha+K}{2\sqrt{2+(\alpha-K)^2}}>0,
     \qquad \alpha\in\RR,
 \]
 we get \ $\lim_{t\uparrow T}C_{\alpha,K}(t)=\infty$, \ and hence
 \[
   \PP\left(\lim_{t\uparrow T}\int_0^t\frac{b(u)^2}{\sigma(u)^2}(X_u^{(\alpha)})^2 \,\dd u=\infty\right)
        =1, \qquad \alpha\in\RR.
 \]
Then by a strong law of large numbers for continuous local martingales, see, e.g.,
 Barczy and Pap \cite[Theorem 15]{BarPap}, we get the MLE of  \ $\alpha$ \ is strongly consistent
 for all \ $\alpha\in\RR$.

The case $K>0$ can be handled in the same way as at the end of the proof of Theorem
 \ref{THM_alpha_1/2_general}. 
\proofend

Finally, we note that in this section we studied the MLE \ $\widehat\alpha_t$ \ of
 \ $\alpha$ \ based on a {\sl continuous} observation \ $(X_s^{(\alpha)})_{s\in[0,t]}$ \
 using the results on Laplace transforms presented in Section \ref{Section_Laplace}.
However, a continuous observation of a diffusion process is only a mathematical idealization,
 in practice the observation is always discrete.
Hence one can pose the question whether our results on the MLE of \ $\alpha$ \ based on
 continuous observations give some information also for discrete observations.
Parameter estimation for discretely observed diffusion processes has been studied by many authors,
 for a detailed discussion and references see, e.g., Bishwal \cite{Bis}.
For discrete observations, one possible approach is to try to find a good
 approximation of the MLE of \ $\alpha$ \
 based on continuous observations (for example, It\^{o} type approximation
 for the stochastic integral in the numerator of \eqref{SEGED20_uj} and
 usual rectangular approximation for the ordinary integral in the denumerator
 of \eqref{SEGED20_uj}).
In this paper we do not consider this question.

\section{$\alpha$-Wiener bridge}\label{Section_Example}

For \ $T\in(0,\infty)$ \ and \ $\alpha\in\RR,$ \ let \ $(X_t^{(\alpha)})_{t\in[0,T)}$ \ be the
 process given by the SDE \eqref{alpha_W_bridge}.
To our knowledge, these kinds of processes in the case of \ $\alpha>0$ \
 have been first considered by Brennan and Schwartz \cite{BreSch}, and see also Mansuy \cite{Man}.
In Brennan and Schwartz \cite{BreSch} the SDE \eqref{alpha_W_bridge} is used to model
 the arbitrage profit associated with a given futures contract in the absence of transaction costs.
By \eqref{SolX}, the unique strong solution of the SDE \eqref{alpha_W_bridge} is
 \begin{align*}
   X_t^{(\alpha)}=\int_{0}^t\left(\frac{T-t}{T-s}\right)^\alpha\,\dd B_s,\qquad t\in[0,T).
 \end{align*}

Theorem \ref{THEOREM_LAPLACE} has the following consequence on the joint Laplace transform of
 \ $\int_0^t \frac{(X_u^{(\alpha)})^2}{(T-u)^2} \, \dd u$ \ and \ $(X_t^{(\alpha)})^2$.

\begin{Thm}\label{Psi}
Let \ $(X_t^{(\alpha)})_{t\in[0,T)}$ \ be the process given by the SDE \eqref{alpha_W_bridge}.
For all \ $\mu > 0$, \ $\nu\geq 0$ \ and \ $t \in [0,T)$, \ we have
 \begin{multline*}
  \EE \exp \left\{ - \mu \int_0^t \frac{(X_u^{(\alpha)})^2}{(T-u)^2} \, \dd u
                   -\nu[X_t^{(\alpha)}]^2
           \right\}
   \\
  = \frac{\left( 1 - \frac{t}{T} \right)^{(1-2\alpha)/4}}
          {\sqrt{\cosh\left(\frac{\sqrt{8 \mu + (2\alpha-1)^2}}{2}
                            \ln\left( 1 - \frac{t}{T} \right)\right)
                 + \frac{1-2\alpha-4\nu(T-t)}{\sqrt{8\mu + (2\alpha-1)^2}}
                   \sinh\left(\frac{\sqrt{8\mu + (2\alpha-1)^2}}{2}
                              \ln\left( 1 - \frac{t}{T} \right)\right)}} .
 \end{multline*}
\end{Thm}

\noindent{\bf Proof.}
Let \ $b(t):=-\frac{1}{T-t}$, \ $t\in[0,T)$, \ and \ $\sigma(t):=1$, \ $t\in[0,T)$.
\ Hence condition \eqref{Laplace_feltetel2} is satisfied with \ $K:=\frac{1}{2}$ \ and
 \ $C:=-\frac{T}{2}$, \ and clearly,
 \[
    B_{K,C}(t)=\left(1-\frac{t}{T}\right)^2,\qquad t\in[0,T).
 \]
By Theorem \ref{THEOREM_LAPLACE}, we have the assertion.
\proofend

Theorem \ref{THM_alpha_1/2_general} has the following consequence on the asymptotic behavior of
 the maximum likelihood estimator \ $\widehat\alpha_t$ \ of \ $\alpha$ \ as \ $t\uparrow T$.

\begin{Thm}\label{THM_maximium_likelihood}
Let \ $(X_t^{(\alpha)})_{t\in[0,T)}$ \ be the process given by the SDE \eqref{alpha_W_bridge}.
For each \ $\alpha>\frac{1}{2}$, \ the maximum likelihood estimator
 \ $\widehat\alpha_t$ \ of \ $\alpha$ \ is asymptotically normal, namely,
 for each \ $\alpha>\frac{1}{2}$,
 \begin{align*}
    \sqrt{I_{\alpha}(t)}(\widehat\alpha_t-\alpha)
      \distr
      \cN(0,1)
      \qquad \text{as \ $t\uparrow T$.}
 \end{align*}
If \ $\alpha = \frac{1}{2},$ \ then the distribution of
 \ $\sqrt{I_{1/2}(t)} \left( \widehat\alpha_t - \frac{1}{2} \right)$ \
 is the same for all \ $t \in (0,T)$, \ namely,
  \[
   \sqrt{I_{1/2}(t)} \left(\widehat\alpha_t - \frac{1}{2} \right)
   \distre
   -\frac{1}{2\sqrt{2}}\,\frac{(W_1)^2-1}{\int_0^1(W_s)^2\,\dd s}
   =
   -\frac{1}{\sqrt{2}} \, \frac{\int_0^1 W_s \, \dd W_s}{\int_0^1 (W_s)^2 \, \dd s},
 \]
 where \ $(W_s)_{s\in[0,1]}$ \ is a standard Wiener process.
\end{Thm}

The following remark is about the asymptotic behavior of the MLE of \ $\alpha$ \ in case
 of \ $\alpha < \frac12$. \ We note that up to our knowledge this case can not be handled
  using only Laplace transforms.

\begin{Rem}\label{Cauchy}
If \ $\alpha < \frac12$, \ then
 \begin{align*}
   \sqrt{I_{\alpha}(t)} \, \big( \widehat\alpha_t - \alpha \big)
   \distr \zeta \qquad \text{as \ $t \uparrow T$,}
 \end{align*}
 where \ $\zeta$ \ is a standard Cauchy distributed random variable, see, e.g.,
 Luschgy \cite[Section 4.2]{Lus1} or Barczy and Pap \cite{BarPap}.
\end{Rem}

Theorem \ref{THM_alpha_1/2_general_random} has the following consequence on
 the (asymptotic) behavior of the MLE of \ $\alpha=1/2$ \ using a {\sl random}
 normalization.

\begin{Thm}\label{THM_maximium_likelihood_rewritten}
Let \ $(X_t^{(\alpha)})_{t\in[0,T)}$ \ be the process given by the SDE \eqref{alpha_W_bridge}.
For all \ $t\in(0,T)$, \ we have
\begin{align*}
   \left(\int_0^t\frac{(X_u^{(1/2)})^2}{(T-u)^2}\,\dd u\right)^{1/2}
      \left(\widehat\alpha_t-\frac{1}{2}\right)
       \distre
         -\frac{\int_0^1 W_s \, \dd W_s}{\left(\int_0^1 (W_s)^2 \, \dd s\right)^{1/2}}
        =-\frac{1}{2}\frac{(W_1)^2-1}{\left(\int_0^1 (W_s)^2 \, \dd s\right)^{1/2}}.
\end{align*}
\end{Thm}

Finally, we note that Es-Sebaiy and Nourdin \cite{EssNou} studied the parameter estimation
 for so-called $\alpha$-fractional bridges which are given by the SDE \eqref{alpha_W_bridge}
 replacing the standard Wiener process \ $B$ \ by a fractional Wiener process.

\par\bigskip\noindent
{\bf Acknowledgment.}
The authors are grateful to the referee for the useful comments.

\end{document}